# WINDING-INVARIANT PRIME IDEALS IN QUANTUM $3 \times 3$ MATRICES

K. R. GOODEARL AND T. H. LENAGAN

ABSTRACT. A complete determination of the prime ideals invariant under winding automorphisms in the generic $3 \times 3$ quantum matrix algebra $\mathcal{O}_q(M_3(k))$ is obtained. Explicit generating sets consisting of quantum minors are given for all of these primes, thus verifying a general conjecture in the $3 \times 3$ case. The result relies heavily on certain tensor product decompositions for winding-invariant prime ideals, developed in an accompanying paper. In addition, new methods are developed here, which show that certain sets of quantum minors, not previously manageable, generate prime ideals in $\mathcal{O}_q(M_n(k))$.

## INTRODUCTION

Although the quantized coordinate ring of $n \times n$ matrices, $\mathcal{O}_q(M_n(k))$, appears to have a very straightforward structure – for instance, it is an iterated skew polynomial ring over the base field – many basic questions about this algebra remain unanswered. In particular, the structure of the prime spectrum of $\mathcal{O}_q(M_n(k))$ is only partially understood, even in the generic case (that is, when $q$ is not a root of unity), where there are far fewer prime ideals than in the root of unity case. The $2 \times 2$ situation is relatively easy and has long since been dealt with; here we complete the picture for the far more complicated $3 \times 3$ quantum matrix algebra.

First, consider the algebra $A = \mathcal{O}_q(M_n(k))$, with arbitrary matrix size $n$ and arbitrary base field $k$, but with $q \in k^\times$ not a root of unity; we recall standard generators and relations for $A$ in (1.1). There is a natural action of the torus $H = (k^\times)^n \times (k^\times)^n$ on $A$ by *winding automorphisms*, which we recall in (1.3). It is known from work of Letzter and the first author that $A$ has only finitely many winding-invariant primes (that is, prime ideals invariant under $H$) [**7**, 5.7(i)], and that the winding-invariant primes control much of the behavior of the prime and primitive spectra of $A$ [**7**, Theorem 6.6]. In particular, spec $A$ has a natural partition into $H$-*strata*, indexed by $H$, such that each $H$-stratum is Zariski-homeomorphic to the prime spectrum of a commutative Laurent polynomial ring over $k$, and such that the primitive ideals within any $H$-stratum are precisely the maximal

This research was partially supported by NSF research grants DMS-9622876 and DMS-9970159 and by NATO Collaborative Research Grant CRG.960250. Some of the research was done while both authors were visiting the Mathematical Sciences Research Institute in Berkeley during the winter of 2000, and they thank MSRI for its support.

Typeset by $\mathcal{A}_{\mathcal{M}}\mathcal{S}$-TEX

elements of that stratum. (See [**1**, Theorems II.2.13, II.8.4, Corollary II.6.5] for a full development of this picture.) Thus, the fundamental step in calculating the prime and primitive spectra of $A$ is to determine the winding-invariant primes. Evidence from the $2 \times 2$ case and from other quantized coordinate rings leads one to the

**Conjecture.** *Every winding-invariant prime of $\mathcal{O}_q(M_n(k))$ can be generated by a set of quantum minors, which can be arranged in a polynormal sequence.*

The main result of this paper is an explicit verification of the conjecture in the case $n = 3$. Cauchon has shown that distinct, comparable winding-invariant primes of $A$ can be separated by the quantum minors they contain [**2**, Proposition 6.2.2 and Théorème 6.3.1], but his methods do not produce sufficiently many quantum minors to provide generating sets. On the other hand, his methods do yield a formula for the number of winding-invariant primes in $A$ (see [**2**, Théorème 3.2.2 and Proposition 3.3.2]); in particular (as we independently verify), there are precisely 230 winding-invariant primes in $\mathcal{O}_q(M_3(k))$. Although our work provides lists of generators (exhibiting a high degree of symmetry) for these primes, we do not yet have either an abstract description or an algorithm for these generating sets. Such a description or algorithm would be essential to further progress – as Cauchon's formula shows, there are 6902 winding-invariant primes in $\mathcal{O}_q(M_4(k))$!

Our procedure may be separated into four distinct parts. The fundamental line, which gets the process started, follows the tensor product decomposition of winding-invariant primes developed in [**6**]; we summarize this line, and establish the requisite notation, in the following section. In particular, this allows us to express every winding-invariant prime of $\mathcal{O}_q(M_3(k))$ as a pullback (via comultiplication and quotient maps) of winding-invariant primes in certain tensor products of 'step-triangular' factor algebras of $\mathcal{O}_q(M_3(k))$. The winding-invariant primes in these factor algebras are built up from winding-invariant primes of $\mathcal{O}_q(M_2(k))$, and so they are easily determined. Inspection of the associated comultiplication and quotient maps then allows us to determine exactly which quantum minors lie in each pullback. Thus, as the second step of our process, we obtain *potential* generating sets for all the winding-invariant primes. What remains is to prove that each of these sets of quantum minors does generate a prime ideal (any ideal generated by quantum minors is automatically winding-invariant), and that the primes thus generated are distinct. In most of these cases, standard skew polynomial technology easily shows that a given set of quantum minors generates a prime ideal. Some cases, however – previewing the difficulties that will need to be overcome in larger quantum matrix algebras – require new methods. Thus, a third part of our work is to establish a new supply of primes generated by quantum minors. This is done with the aid of our fourth line, which develops three new bases for $\mathcal{O}_q(M_n(k))$, similar to the *standard basis* of [**11**] and the *preferred basis* of [**5**]. (That portion of our work, in Section 5, applies as well to the root of unity case as to the generic case.)

## 1. Winding-invariant primes and tensor product decompositions

We begin by establishing our notation and outlining the tools from [**6**] that we shall use. Since most of our work concerns quantized coordinate rings of square matrices, we restrict attention to that case in this section. On the occasions when rectangular quantum matrix algebras are needed, we will use analogous notation.

**1.1.** Fix a base field $k$, a nonzero scalar $q \in k^\times$ which is not a root of unity, and a positive integer $n$. Set $A = \mathcal{O}_q(M_n(k))$, presented by generators $X_{ij}$ for $i, j = 1, \ldots, n$ and relations

$$\begin{aligned}
X_{ij}X_{lj} &= qX_{lj}X_{ij} & (i < l) \\
X_{ij}X_{im} &= qX_{im}X_{ij} & (j < m) \\
X_{ij}X_{lm} &= X_{lm}X_{ij} & (i < l, \, j > m) \\
X_{ij}X_{lm} - X_{lm}X_{ij} &= (q - q^{-1})X_{im}X_{lj} & (i < l, \, j < m).
\end{aligned}$$

Recall that $A$ is a bialgebra with comultiplication $\Delta : A \to A \otimes A$ and counit $\epsilon : A \to k$ such that

$$\Delta(X_{ij}) = \sum_{l=1}^{n} X_{il} \otimes X_{lj} \qquad \text{and} \qquad \epsilon(X_{ij}) = \delta_{ij}$$

for all $i, j$.

It is well known that $A$ can be presented as an iterated skew polynomial algebra over $k$, with the generators adjoined in lexicographic order. (Cf. [**8**, pp. 10, 11] for notation and conventions concerning skew polynomial rings.) Thus,

$$A = k[X_{11}][X_{12}; \tau_{12}] \cdots [X_{ij}; \tau_{ij}, \delta_{ij}] \cdots [X_{nn}; \tau_{nn}, \delta_{nn}]$$

where each $\tau_{lm}$ is a $k$-algebra automorphism and each $\delta_{lm}$ is a $k$-linear $\tau_{lm}$-derivation, such that

$$\tau_{lm}(X_{ij}) = \begin{cases} q^{-1}X_{ij} & (i = l, \, j < m) \\ q^{-1}X_{ij} & (i < l, \, j = m) \\ X_{ij} & (i < l, \, j \neq m) \end{cases}$$

$$\delta_{lm}(X_{ij}) = \begin{cases} 0 & (i = l, \, j < m) \\ 0 & (i < l, \, j \geq m) \\ (q^{-1} - q)X_{im}X_{lj} & (i < l, \, j < m). \end{cases}$$

**1.2.** As in [**5**], we write quantum minors in the form $[X \mid Y]$ where $X$ and $Y$ are either sets or lists of row and column indices, but here we typically write both row and column indices in ascending order. More precisely, let $I = \{i_1, \ldots, i_t\}$ and $J = \{j_1, \ldots, j_t\}$ be $t$-element subsets of $\{1, \ldots, n\}$, arranged so that $i_1 < i_2 < \cdots < i_t$ and $j_1 < \cdots < j_t$. There is a natural $k$-algebra embedding $\mathcal{O}_q(M_t(k)) \to A$ such that $X_{rc} \mapsto X_{i_r j_c}$ for $r, c = 1, \ldots, t$, and the image of the quantum determinant of $\mathcal{O}_q(M_t(k))$ under this embedding is the $t \times t$ quantum minor labelled $[I \mid J]$ or $[i_1 \cdots i_t \mid j_1 \cdots j_t]$. Recall that comultiplication of quantum minors is given by the rule

$$\Delta\big([I \mid J]\big) = \sum_{\substack{K \subseteq \{1,\ldots,n\} \\ |K| = |I|}} [I \mid K] \otimes [K \mid J]$$

[**14**, Equation 1.9].

**1.3.** Set $H = (k^\times)^n \times (k^\times)^n$, and recall that $H$ acts on $A$ by $k$-algebra automorphisms such that
$$(\alpha_1, \ldots, \alpha_n, \beta_1, \ldots, \beta_n).X_{ij} = \alpha_i \beta_j X_{ij}$$
for $(\alpha_1, \ldots, \alpha_n, \beta_1, \ldots, \beta_n) \in H$ and all $i, j$. Although this action is not faithful, it allows a convenient notation. The automorphisms affording this action arise from winding automorphisms of $A$, as follows.

If $\phi$ is any character of $A$, that is, a $k$-algebra homomorphism $A \to k$, the rules
$$\tau_\phi^l(a) = \sum_{(a)} \phi(a_1) a_2 \qquad\qquad \tau_\phi^r(a) = \sum_{(a)} a_1 \phi(a_2)$$
define $k$-algebra endomorphisms $\tau_\phi^l$ and $\tau_\phi^r$ of $A$. When $\phi$ is invertible in $A^*$ (with respect to the convolution product), $\tau_\phi^l$ and $\tau_\phi^r$ are automorphisms of $A$, called the *left and right winding automorphisms* arising from $\phi$. As is easily checked, if $\phi$ is convolution-invertible then $\phi(X_{ij}) = 0$ for $i \neq j$, and the scalars $\phi_i = \phi(X_{ii})$ must all be nonzero. Hence, $\tau_\phi^l$ and $\tau_\phi^r$ are given by the actions of the elements
$$(\phi_1, \ldots, \phi_n, 1, \ldots, 1) \qquad \text{and} \qquad (1, \ldots, 1, \phi_1, \ldots, \phi_n)$$
from $H$.

In particular, the ideals of $A$ invariant under all winding automorphisms are just the $H$-invariant ideals, and so we restate our goal as that of determining the $H$-invariant prime ideals of $A$. Recall that an *$H$-prime ideal* of $A$ is any proper $H$-invariant ideal $P$ such that $P$ contains a product $IJ$ of $H$-invariant ideals $I$ and $J$ only when $P$ contains $I$ or $J$. Obviously $H$-invariant primes are $H$-prime; the converse holds in $A$ by [**7**, (5.7)(i)], where it is also proved that $A$ has only finitely many $H$-primes (in fact, at most $2^{n^2}$ of them).

**1.4.** Let $\boldsymbol{RC}$ denote the set of all pairs $(\boldsymbol{r}, \boldsymbol{c})$ where $\boldsymbol{r} = (r_1, \ldots, r_t)$ and $\boldsymbol{c} = (c_1, \ldots, c_t)$ are strictly increasing sequences in $\{1, \ldots, n\}$ of the same length. We allow $t = 0$, in which case $\boldsymbol{r}$ and $\boldsymbol{c}$ are empty sequences. Write $(\boldsymbol{r}, \boldsymbol{c}) \in \boldsymbol{RC}_t$ to indicate that $\boldsymbol{r}$ and $\boldsymbol{c}$ have length $t$. We shall partition the prime spectrum of $A$ into subsets indexed by $\boldsymbol{RC}$. Part of the labelling for this partition requires the following partial order on index sets.

Let $I, I' \subseteq \{1, \ldots, n\}$ be index sets of the same cardinality, with elements listed in ascending order, say $I = \{i_1 < \cdots < i_l\}$ and $I' = \{i'_1 < \cdots < i'_l\}$ for short. Define $I \leq I'$ to mean that $i_s \leq i'_s$ for $s = 1, \ldots, l$. (In Section 5, where several order relations on index sets come into play, the one just described will be denoted $\leq_c$.)

For $(\boldsymbol{r}, \boldsymbol{c}) \in \boldsymbol{RC}_t$, let $K_{\boldsymbol{r}, \boldsymbol{c}}$ be the ideal of $A$ generated by the following set of quantum minors:
$$\{[I \mid J] \mid |I| > t\} \cup \{[I \mid J] \mid |I| = l \leq t \text{ and } I \not\geq \{r_1, \ldots, r_l\}\}$$
$$\cup \{[I \mid J] \mid |I| = l \leq t \text{ and } J \not\geq \{c_1, \ldots, c_l\}\}.$$

Now set $d_l^{\boldsymbol{r}, \boldsymbol{c}} = [r_1 \cdots r_l \mid c_1 \cdots c_l]$ for $l \leq t$, and let $D_{\boldsymbol{r}, \boldsymbol{c}}$ denote the multiplicative subset of $A$ generated by $d_1^{\boldsymbol{r}, \boldsymbol{c}}, \ldots, d_t^{\boldsymbol{r}, \boldsymbol{c}}$. Finally, set
$$\operatorname{spec}_{\boldsymbol{r}, \boldsymbol{c}} A = \{P \in \operatorname{spec} A \mid K_{\boldsymbol{r}, \boldsymbol{c}} \subseteq P \text{ and } P \cap D_{\boldsymbol{r}, \boldsymbol{c}} = \varnothing\}$$

and $H\text{-}\operatorname{spec}_{\boldsymbol{r},\boldsymbol{c}} A = (H\text{-}\operatorname{spec} A) \cap (\operatorname{spec}_{\boldsymbol{r},\boldsymbol{c}} A)$, where $\operatorname{spec} A$ and $H\text{-}\operatorname{spec} A$ denote the collections of prime and $H$-prime ideals of $A$, respectively.

**Theorem.** *The set* $\operatorname{spec} A$ *is the disjoint union of the subsets* $\operatorname{spec}_{\boldsymbol{r},\boldsymbol{c}} A$ *for* $(\boldsymbol{r}, \boldsymbol{c}) \in \boldsymbol{RC}$. *Consequently,* $H\text{-}\operatorname{spec} A$ *is the disjoint union of its subsets* $H\text{-}\operatorname{spec}_{\boldsymbol{r},\boldsymbol{c}} A$.

*Proof.* The first statement is [**6**, Corollary 1.10]; the second follows from that together with the fact (noted above) that all $H$-primes of $A$ are prime. □

**1.5.** Given $(\boldsymbol{r}, \boldsymbol{c}) \in \boldsymbol{RC}_t$ for some $t$, set

$$R^+_{\boldsymbol{r},0} = A/\langle X_{ij} \mid j > t \text{ or } i < r_j \rangle \quad \text{and} \quad R^-_{\boldsymbol{c},0} = A/\langle X_{ij} \mid i > t \text{ or } j < c_i \rangle.$$

Write $Y_{ij}$ and $Z_{ij}$ for the images of $X_{ij}$ in $R^+_{\boldsymbol{r},0}$ and $R^-_{\boldsymbol{c},0}$, respectively. Note that these algebras are iterated skew polynomial extensions of $k$, hence noetherian domains, the natural indeterminates for these iterated skew polynomial structures being those $Y_{ij}$ and $Z_{ij}$ which are nonzero. The $Y_{r_s s}$ are regular normal elements in $R^+_{\boldsymbol{r},0}$, and the $Z_{sc_s}$ are regular normal elements in $R^-_{\boldsymbol{c},0}$ (cf. [**6**, (2.1)]). Hence, we can form Ore localizations

$$R^+_{\boldsymbol{r}} = R^+_{\boldsymbol{r},0}[Y^{-1}_{r_1 1}, \ldots, Y^{-1}_{r_t t}] \quad \text{and} \quad R^-_{\boldsymbol{c}} = R^-_{\boldsymbol{c},0}[Z^{-1}_{1c_1}, \ldots, Z^{-1}_{tc_t}].$$

These algebras are mixed iterated skew polynomial and skew-Laurent extensions of $k$, hence noetherian domains.

**1.6.** Given $(\boldsymbol{r}, \boldsymbol{c}) \in \boldsymbol{RC}$, let $\pi^+_{\boldsymbol{r},0} : A \to R^+_{\boldsymbol{r},0}$ and $\pi^-_{\boldsymbol{c},0} : A \to R^-_{\boldsymbol{c},0}$ be the quotient maps, and define

$$\beta_{\boldsymbol{r},\boldsymbol{c}} : A \xrightarrow{\Delta} A \otimes A \xrightarrow{\pi^+_{\boldsymbol{r},0} \otimes \pi^-_{\boldsymbol{c},0}} R^+_{\boldsymbol{r},0} \otimes R^-_{\boldsymbol{c},0} \xrightarrow{\subseteq} R^+_{\boldsymbol{r}} \otimes R^-_{\boldsymbol{c}}.$$

Thus $\beta_{\boldsymbol{r},\boldsymbol{c}}$ is a $k$-algebra homomorphism, and it satisfies

$$\beta_{\boldsymbol{r},\boldsymbol{c}}(X_{ij}) = \sum_{l \leq t,\, r_l \leq i,\, c_l \leq j} Y_{il} \otimes Z_{lj}$$

for all $i, j$. In particular, $\beta_{\boldsymbol{r},\boldsymbol{c}}(X_{ij}) = 0$ when $i < r_1$ or $j < c_1$.

**Theorem.** *For each* $(\boldsymbol{r}, \boldsymbol{c}) \in \boldsymbol{RC}$, *there is a bijection*

$$(H\text{-}\operatorname{spec} R^+_{\boldsymbol{r}}) \times (H\text{-}\operatorname{spec} R^-_{\boldsymbol{c}}) \longrightarrow H\text{-}\operatorname{spec}_{\boldsymbol{r},\boldsymbol{c}} A$$

*given by the rule* $(Q^+, Q^-) \mapsto \beta^{-1}_{\boldsymbol{r},\boldsymbol{c}}\big((Q^+ \otimes R^-_{\boldsymbol{c}}) + (R^+_{\boldsymbol{r}} \otimes Q^-)\big)$.

*Proof.* [**6**, Theorem 3.5]. □

For $Q^\pm$ as in the theorem, the $H$-prime $P = \beta^{-1}_{\boldsymbol{r},\boldsymbol{c}}\big((Q^+ \otimes R^-_{\boldsymbol{c}}) + (R^+_{\boldsymbol{r}} \otimes Q^-)\big)$ equals the kernel of the map

$$A \xrightarrow{\beta_{\boldsymbol{r},\boldsymbol{c}}} R^+_{\boldsymbol{r}} \otimes R^-_{\boldsymbol{c}} \xrightarrow{\pi_{Q^+} \otimes \pi_{Q^-}} (R^+_{\boldsymbol{r}}/Q^+) \otimes (R^-_{\boldsymbol{c}}/Q^-),$$

where $\pi_{Q^+}$ and $\pi_{Q^-}$ are quotient maps. Thus, we refer to this description of $P$ as a *tensor product decomposition*.

2. SETTING UP THE SEARCH FOR $H$-PRIMES IN $\mathcal{O}_q(M_3(k))$

In case $n = 3$, the algebras $R_{\boldsymbol{r}}^+$ and $R_{\boldsymbol{c}}^-$ defined in (1.5) are small and easily analyzed; in particular, their $H$-primes can be readily determined. That allows us, via Theorem 1.6, to write down an explicit list of the $H$-primes of $\mathcal{O}_q(M_3(k))$. However, in such a list the $H$-primes appear in the form $\beta_{\boldsymbol{r},\boldsymbol{c}}^{-1}\big((Q^+ \otimes R_{\boldsymbol{c}}^-) + (R_{\boldsymbol{r}}^+ \otimes Q^-)\big)$, and we do not yet have the technology to directly compute generating sets for these ideals in all cases. Thus, we partly take a roundabout route, in which we find sufficiently many $H$-primes (with explicit generating sets) "similar to" the above ideals to ensure that all the $H$-primes of $\mathcal{O}_q(M_3(k))$ are accounted for. In this section, we take the first steps along the above route, establishing some convenient notation and recording the easiest cases. In these cases, where $\boldsymbol{r}$ and $\boldsymbol{c}$ have length 0 or 1, the $H$-primes of $\mathcal{O}_q(M_3(k))$ are already known, along with convenient sets of generators.

**2.1.** Throughout this section, fix $n = 3$ and keep the notation of Section 1. The possible sequences $\boldsymbol{r}$ and $\boldsymbol{c}$ with which to build pairs $(\boldsymbol{r}, \boldsymbol{c}) \in \boldsymbol{RC}$ are as follows:

$$(1,2,3),\ (1,2),\ (1,3),\ (2,3),\ (1),\ (2),\ (3),\ \varnothing.$$

The corresponding eight possibilities for each of $R_{\boldsymbol{r}}^+$ and $R_{\boldsymbol{c}}^-$ are listed in Figure 2.1 below, with abbreviated notation. In the given $3 \times 3$ patterns, the symbols $0$, $+$, $\pm$ indicate that the coset $Y_{ij}$ or $Z_{ij}$ corresponding to that position is zero, or nonzero, or nonzero and inverted, respectively. Since $Y_{ij} = 0$ for $i < j$ in all cases, we omit the zeroes above the diagonal in the descriptions of the $R_{\boldsymbol{r}}^+$, and similarly below the diagonal for the $R_{\boldsymbol{c}}^-$.

$$R_{123}^+ = k \begin{bmatrix} \pm & & \\ + & \pm & \\ + & + & \pm \end{bmatrix} \qquad R_{12}^+ = k \begin{bmatrix} \pm & & \\ + & \pm & \\ + & + & 0 \end{bmatrix} \qquad R_{13}^+ = k \begin{bmatrix} \pm & & \\ + & 0 & \\ + & \pm & 0 \end{bmatrix} \qquad R_{23}^+ = k \begin{bmatrix} 0 & & \\ \pm & 0 & \\ + & \pm & 0 \end{bmatrix}$$

$$R_1^+ = k \begin{bmatrix} \pm & & \\ + & 0 & \\ + & 0 & 0 \end{bmatrix} \qquad R_2^+ = k \begin{bmatrix} 0 & & \\ \pm & 0 & \\ + & 0 & 0 \end{bmatrix} \qquad R_3^+ = k \begin{bmatrix} 0 & & \\ 0 & 0 & \\ \pm & 0 & 0 \end{bmatrix} \qquad R_\varnothing^+ = k \begin{bmatrix} 0 & & \\ 0 & 0 & \\ 0 & 0 & 0 \end{bmatrix}$$

$$R_{123}^- = k \begin{bmatrix} \pm & + & + \\ & \pm & + \\ & & \pm \end{bmatrix} \qquad R_{12}^- = k \begin{bmatrix} \pm & + & + \\ & \pm & + \\ & & 0 \end{bmatrix} \qquad R_{13}^- = k \begin{bmatrix} \pm & + & + \\ & 0 & \pm \\ & & 0 \end{bmatrix} \qquad R_{23}^- = k \begin{bmatrix} 0 & \pm & + \\ & 0 & \pm \\ & & 0 \end{bmatrix}$$

$$R_1^- = k \begin{bmatrix} \pm & + & + \\ & 0 & 0 \\ & & 0 \end{bmatrix} \qquad R_2^- = k \begin{bmatrix} 0 & \pm & + \\ & 0 & 0 \\ & & 0 \end{bmatrix} \qquad R_3^- = k \begin{bmatrix} 0 & 0 & \pm \\ & 0 & 0 \\ & & 0 \end{bmatrix} \qquad R_\varnothing^- = k \begin{bmatrix} 0 & 0 & 0 \\ & 0 & 0 \\ & & 0 \end{bmatrix}$$

Figure 2.1

**2.2.** Note that each of the algebras $R_\bullet^\pm$ is a (possibly iterated) skew-Laurent extension of a localized factor of $\mathcal{O}_q(M_2(k))$; for instance,

$$R_{123}^- \cong \big(\mathcal{O}_q(M_2(k))[X_{21}^{-1}]\big)[z_1^{\pm 1}, z_2^{\pm 1}; \sigma_1, \sigma_2]$$

for two (commuting) automorphisms $\sigma_i$ of $\mathcal{O}_q(M_2(k))[X_{21}^{-1}]$. (See, e.g., [**8**, pp. 16, 17] for notation and details concerning skew-Laurent extensions.) The following lemma shows that the $H$-primes of $R_{123}^-$ are all induced from $H$-primes of $\mathcal{O}_q(M_2(k))[X_{21}^{-1}]$.

**Lemma.** *Let $T = S[z^{\pm 1}; \sigma]$ where $S$ is a $k$-algebra and $\sigma$ is a $k$-algebra automorphism of $S$. Let $G$ be a group acting on $T$ by $k$-algebra automorphisms, such that $S$ is $G$-invariant and $\sigma$ coincides with the action of some element of $G$ on $S$. Assume also that there exists $g_0 \in G$ such that $g_0$ acts trivially on $S$ while $g_0(z) = rz$ for some non-root of unity $r \in k^\times$. Then the $G$-primes of $T$ are exactly the ideals induced from $G$-primes of $S$.*

*Proof.* Note that if $I$ is a $G$-invariant ideal of $S$, then $I$ is also $\sigma$-invariant, and so $IT = TI$ is a $G$-invariant ideal of $T$. We claim that all $G$-invariant ideals of $T$ have this form. Thus, let $P$ be an arbitrary $G$-invariant ideal of $T$, and set $I = P \cap S$, a $G$-invariant ideal of $S$. Since $T/IT \cong (S/I)[z^{\pm 1}; \sigma]$, in proving the claim it suffices to assume that $P \cap S = 0$ and show that $P = 0$.

If $P \ne 0$, choose a nonzero element $p \in P$ of minimal length, say length $n+1$. After multiplying $p$ by a suitable power of $z$, we may assume that $p = p_0 + p_1 z + \cdots + p_n z^n$ for some $p_i \in S$, where $p_0, p_n \ne 0$. Since $P$ is $G$-invariant, it also contains the polynomial $g_0(p) = p_0 + rp_1 z + \cdots + r^n p_n z^n$. The difference $g_0(p) - p$ is then an element of $P$ of length at most $n$, whence $g_0(p) - p = 0$. Now $r^n p_n = p_n$. Since $r$ is not a root of unity and $p_n \ne 0$, we must have $n = 0$. But then $p = p_0 \in P \cap S$, contradicting our assumption that $P \cap S = 0$. Therefore $P = 0$, establishing the claim.

Now let $P$ be a $G$-prime of $T$. By the claim, $P = QT$ for some proper $G$-invariant ideal $Q$ of $S$. If $I$ and $J$ are $G$-invariant ideals of $S$ with $IJ \subseteq Q$, then $IT = TI$ and $JT = TJ$ are $G$-invariant ideals of $T$ with $(IT)(JT) = IJT \subseteq P$, whence $P$ contains $IT$ or $JT$, and so $Q$ contains $I$ or $J$. Thus $Q$ is a $G$-prime of $S$.

Conversely, let $Q$ be any $G$-prime of $S$. Then $QT$ is a proper $G$-invariant ideal of $T$. Suppose that $I'$ and $J'$ are $G$-invariant ideals of $T$ such that $I'J' \subseteq QT$. Then $I' = IT = TI$ and $J' = JT = TJ$ for some $G$-invariant ideals $I$ and $J$ of $S$. Further, $IJT = I'J' \subseteq QT$ and so $IJ \subseteq Q$, whence $Q$ contains $I$ or $J$, and thus $QT$ contains $I'$ or $J'$. Therefore $QT$ is a $G$-prime of $T$. □

**2.3.** The $H$-primes in $\mathcal{O}_q(M_2(k))$ have long been known; they can, for example, be determined very quickly from Theorem 1.6, as noted in [**6**, (4.1)]. There are 14 of these $H$-primes, and we can give generating sets for them in abbreviated form as in the following display:

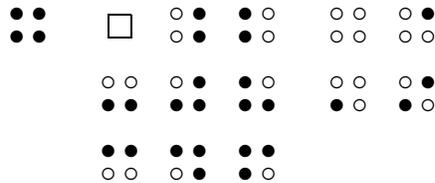

Here the circles are placeholders, a bullet in position $(i,j)$ is a marker for a generator $X_{ij}$, and the square denotes the $2 \times 2$ quantum determinant in $\mathcal{O}_q(M_2(k))$. For example, the schematic ○●/○○ indicates the ideal $\langle X_{12}, X_{22}\rangle$, while ○○/○○ denotes the zero ideal.

With the aid of Lemma 2.2, we can now write down all the $H$-primes in the algebras $R_r^+$ and $R_c^-$. For instance, the $H$-primes of $R_{123}^-$ are induced from $H$-primes of $\mathcal{O}_q(M_2(k))[X_{21}^{-1}]$, and thus from $H$-primes of $\mathcal{O}_q(M_2(k))$ not containing $X_{21}$; there are 6 such $H$-primes. The $H$-primes in each $R_c^-$ are symbolized in Figure 2.3 below; those for the $R_r^+$ are obtained by transposition. The symbols $\circ$, $\bullet$, and $\square$ are used as above, except that now the square denotes the $2 \times 2$ quantum minor $[12 \mid 23]$. We use asterisks for placeholders in the positions $(i, c_i)$ as reminders that the generators $Z_{ic_i}$ are inverted in $R_c^-$ and thus are not candidates for generators of $H$-primes.

$$
\begin{array}{cl}
H\text{-spec } R_{123}^- : & *\circ\circ \quad *\square \quad *\circ\bullet \quad *\bullet\bullet \quad *\circ\bullet \quad *\bullet\bullet \\
 & *\circ \quad *\square \quad *\circ \quad *\circ \quad *\bullet \quad *\bullet \\
 & * \quad * \quad * \quad * \quad * \quad * \\[4pt]
H\text{-spec } R_{12}^- : & *\circ\circ \quad *\square \quad *\circ\bullet \quad *\bullet\bullet \quad *\circ\bullet \quad *\bullet\bullet \\
 & *\circ \quad *\square \quad *\circ \quad *\circ \quad *\bullet \quad *\bullet \\[4pt]
H\text{-spec } R_{13}^- : & *\circ\circ \quad *\circ\bullet \quad *\bullet\circ \quad *\bullet\bullet \\
 & * \quad * \quad * \quad * \\[4pt]
H\text{-spec } R_{23}^- : & *\circ \quad *\bullet \\
 & * \quad * \\[4pt]
H\text{-spec } R_1^- : & *\circ\circ \quad *\circ\bullet \quad *\bullet\circ \quad *\bullet\bullet \\
H\text{-spec } R_2^- : & *\circ \quad *\bullet \\
H\text{-spec } R_3^- : & * \\
H\text{-spec } R_\varnothing^- : & \circ
\end{array}
$$

Figure 2.3

**2.4.** When $r = c = \varnothing$, the only $H$-primes in $R_r^+$ and $R_c^-$ are the zero ideals, and $\beta_{r,c}^{-1}(\langle 0 \rangle) = \langle X_{ij} \mid i, j = 1, 2, 3 \rangle$. We record this $H$-prime in Figure 2.4.

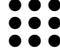

Figure 2.4

**2.5.** Next, consider the cases when $r, c = (1), (2), (3)$. In these cases, Theorem 1.6 yields those $H$-primes of $A$ which contain all the $2 \times 2$ quantum minors but not all the $X_{ij}$. These $H$-primes were determined in [**4**]; there are 49, each generated by the $2 \times 2$ quantum minors together with some of the sets $\{X_{i1}, X_{i2}, X_{i3}\}$ and/or $\{X_{1j}, X_{2j}, X_{3j}\}$. If we express these $H$-primes in the form $\beta_{r,c}^{-1}\big((Q^+ \otimes R_c^-) + (R_r^+ \otimes Q^-)\big)$, generating sets may be displayed as in Figure 2.5 below, where the rows and columns are indexed by the $Q^+$ and $Q^-$, respectively. (The schematics for the $Q^-$ are taken from Figure 2.3, and those for the $Q^+$ are obtained by transposition.) Here the squares and rectangles stand for the corresponding $2 \times 2$ quantum minors; we have not drawn those minors which already lie in the ideals generated by the marked $X_{ij}$. For instance, the schematic in the third position of the second row indicates the ideal $\langle X_{12}, X_{22}, X_{31}, X_{32}, X_{33}, [12|13] \rangle$; note that all $2 \times 2$ quantum minors do lie in this ideal. Two squares with a common edge indicate <u>three</u> quantum minors, symbolized by two squares and one rectangle. For instance, the second schematic in the first row stands for the ideal $\langle [12|12], [13|12], [23|12], X_{13}, X_{23}, X_{33} \rangle$.

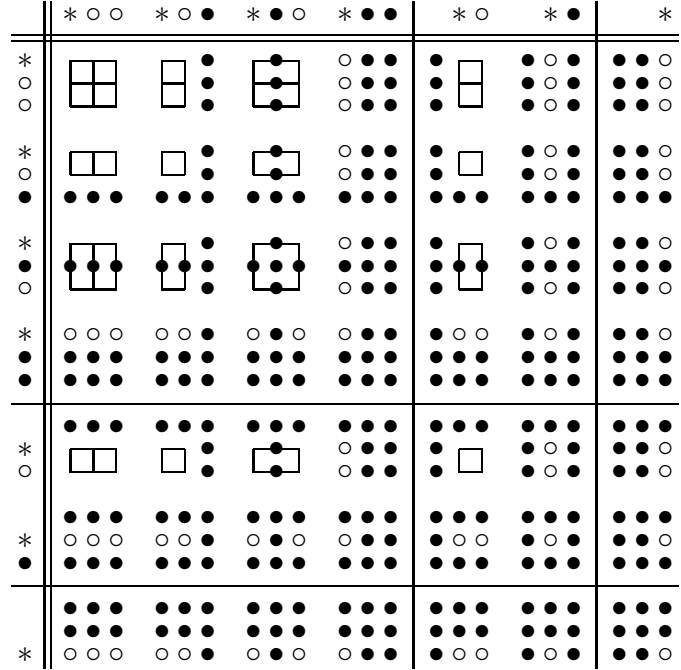

Figure 2.5

## 3. $H$-PRIMES NOT CONTAINING THE QUANTUM DETERMINANT

We continue the analysis of $A = \mathcal{O}_q(M_3(k))$ begun in the previous section.

**3.1.** Let us next consider the $H$-primes of $A$ which do not contain the quantum determinant; these are the $H$-primes which lie in $\operatorname{spec}_{\boldsymbol{r},\boldsymbol{c}} A$ for $\boldsymbol{r} = \boldsymbol{c} = (1,2,3)$. As indicated in Figure 2.3, each of $R_{123}^{\pm}$ has exactly six $H$-primes. Hence, Theorem 1.6 shows that there are exactly 36 $H$-primes in $A$ not containing the quantum determinant. In fact, we can verify this count using earlier results, as follows. First, by standard localization techniques, the $H$-primes in $A$ not containing the quantum determinant are in bijection with the $H$-primes of $\mathcal{O}_q(GL_3(k))$ (cf. [**1**, Exercise II.1.J]). Second, it is known that the $H$-primes of $\mathcal{O}_q(GL_3(k))$ are in bijection with the $H_1$-primes of $\mathcal{O}_q(SL_3(k))$ for a suitable torus $H_1$ (see [**1**, Lemma II.5.16]). Finally, the $H_1$-primes of $\mathcal{O}_q(SL_3(k))$ can be computed by the methods of Hodges and Levasseur [**9**]; there are exactly 36 of them, indexed by $S_3 \times S_3$. (Alternatively, one can write $\mathcal{O}_q(GL_3(k))$ as a factor algebra of $\mathcal{O}_q(SL_4(k))$ and compute the corresponding $H_2$-primes of $\mathcal{O}_q(SL_4(k))$, for an appropriate torus $H_2$, by the methods of [**10**].) Such an analysis determines the $H$-primes of $\mathcal{O}_q(M_3(k))$ not containing the quantum determinant as contractions of $H$-primes of $\mathcal{O}_q(GL_3(k))$, and so only determines generators for these $H$-primes up to torsion with respect to powers of the quantum determinant. Below, we find actual generating sets.

Since we do not have appropriate methods to easily compute these $H$-primes in the form given by Theorem 1.6, we proceed as follows:

(1) For each pair $(Q^+, Q^-)$ with $Q^{\pm} \in H\text{-}\operatorname{spec} R_{123}^{\pm}$, we compute the quantum minors

lying in the ideal $\beta_{r,c}^{-1}((Q^+\otimes R_c^-)+(R_r^+\otimes Q^-))$, and we consider that set of quantum minors as a generating set for an (a priori smaller) ideal. These generating sets are displayed in Figure 3.1 below.
  (2) We prove that each generating set in Figure 3.1 generates an $H$-prime of $A$.
  (3) We check that the $H$-primes generated in Step 2 are all distinct.

Once Step 3 is complete, we will have produced 36 distinct $H$-primes in $\operatorname{spec}_{r,c} A$, and since this set has only 36 members, we will have determined it.

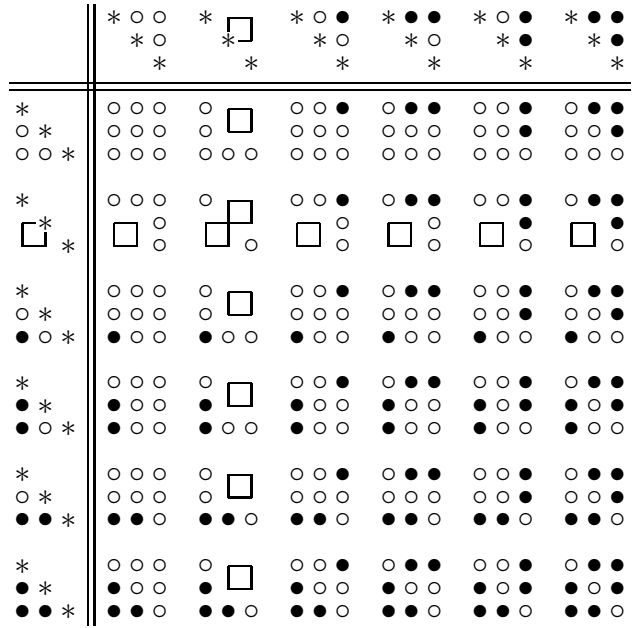

Figure 3.1

Obviously all the ideals corresponding to Figure 3.1 are $H$-invariant. Many of them can be shown to be completely prime via the following general lemma.

**3.2. Lemma.** *Let $T = \mathcal{O}_q(M_{m,n}(k))$ and $M \times N = \{1,\ldots,m\} \times \{1,\ldots,n\}$ for some positive integers $m$, $n$, and let $\mathcal{C}$ be a subset of $M \times N$ such that*

  (*) *For $1 \le i < s \le m$ and $1 \le j < t \le n$, if $(i,j)$ or $(s,t)$ is in $\mathcal{C}$, then $(i,t)$ or $(s,j)$ is in $\mathcal{C}$.*

*Set $P = \langle X_{ij} \mid (i,j) \in \mathcal{C}\rangle$, and set $Y_{ij} = X_{ij} + P$ for all $i,j$. Then $T/P$ is an iterated skew polynomial algebra over $k$ with respect to the variables $Y_{ij}$ for $(i,j) \in (M \times N) \setminus \mathcal{C}$ (ordered lexicographically). In particular, $P$ is a completely prime ideal of $T$.*

*Proof.* Give $M \times N$ the lexicographic order, and view $T$ in the standard way as an iterated skew polynomial algebra over $k$ with respect to the variables $X_{ij}$ (cf. (1.1) for the case where $m = n$). For $(s,t) \in M \times N$, set

$$R_{st} = k\langle X_{ij} \mid (i,j) <_{\text{lex}} (s,t)\rangle$$
$$Q_{st} = \langle X_{ij} \mid (i,j) \in \mathcal{C} \text{ and } (i,j) <_{\text{lex}} (s,t)\rangle \triangleleft R_{st}$$
$$T_{st} = k\langle X_{ij} \mid (i,j) \leq_{\text{lex}} (s,t)\rangle$$
$$P_{st} = \langle X_{ij} \mid (i,j) \in \mathcal{C} \text{ and } (i,j) \leq_{\text{lex}} (s,t)\rangle \triangleleft T_{st}.$$

Thus $R_{11} = k$ and $T_{mn} = T$, and each $T_{st} = R_{st}[X_{st}; \tau_{st}, \delta_{st}]$ for some automorphism $\tau_{st}$ and $\tau_{st}$-derivation $\delta_{st}$. We claim that

(1) $Q_{st}$ is a $(\tau_{st}, \delta_{st})$-ideal of $R_{st}$;
(2) If $(s,t) \notin \mathcal{C}$, then $P_{st} = Q_{st}T_{st}$, while if $(s,t) \in \mathcal{C}$, then $P_{st} = Q_{st}T_{st} + T_{st}X_{st}$;
(3) $P_{st} \cap R_{st} = Q_{st}$

for all $(s,t) \in M \times N$.

Note first that all $X_{ij}$ for $(i,j) <_{\text{lex}} (s,t)$ are $\tau_{st}$-eigenvectors, whence $Q_{st}$ is $\tau_{st}$-stable. To prove $\delta_{st}$-stability, it suffices to show that $\delta_{st}(X_{ij}) \in Q_{st}$ for $(i,j) \in \mathcal{C}$ with $(i,j) <_{\text{lex}} (s,t)$. This is clear if $i = s$ or $j \geq t$, since then $\delta_{st}(X_{ij}) = 0$. If $i < s$ and $j < t$, then $\delta_{st}(X_{ij}) = (q^{-1} - q)X_{it}X_{sj}$. Hypothesis (*) implies that $X_{it} \in Q_{st}$ or $X_{sj} \in Q_{st}$, whence $\delta_{st}(X_{ij}) \in Q_{st}$ as desired. This proves claim (1).

In view of (1), $Q_{st}T_{st}$ is an ideal of $T_{st}$. Thus, the first case of (2) is clear. Now assume that $(s,t) \in \mathcal{C}$; it suffices to show that $Q_{st}T_{st} + T_{st}X_{st}$ is an ideal of $T_{st}$. Hence, it is enough to show that $X_{st}X_{ij} \in Q_{st}T_{st} + T_{st}X_{st}$ for all $(i,j) <_{\text{lex}} (s,t)$. This is clear if $i = s$ or $j \geq t$, since then $X_{st}X_{ij}$ is a scalar multiple of $X_{ij}X_{st}$. If $i < s$ and $j < t$, then $X_{st}X_{ij} = X_{ij}X_{st} + (q^{-1} - q)X_{it}X_{sj}$. By (*), $X_{it} \in Q_{st}$ or $X_{sj} \in Q_{st}$, and again $X_{st}X_{ij} \in Q_{st}T_{st} + T_{st}X_{st}$. Thus claim (2) is proved. Claim (3) follows immediately.

By induction on (3), we see that $P \cap R_{st} = Q_{st}$ and $P \cap T_{st} = P_{st}$ for all $(s,t) \in M \times N$. Hence, there are natural $k$-algebra isomorphisms

$$R_{st}/Q_{st} \cong \overline{R}_{st} = (R_{st} + P)/P \qquad \text{and} \qquad T_{st}/P_{st} \cong \overline{T}_{st} = (T_{st} + P)/P.$$

When $(s,t) \in \mathcal{C}$, we have $X_{st} \in P$ and so $\overline{T}_{st} = \overline{R}_{st}$. On the other hand, when $(s,t) \notin \mathcal{C}$ it follows from (1) and (2) that $\overline{T}_{st} = \overline{R}_{st}[\overline{X}_{st}; \tau_{st}, \delta_{st}]$, a skew polynomial extension. The lemma follows. $\square$

**3.3.** Lemma 3.2 immediately shows that 25 of the displayed schematics in Figure 3.1 – those in which no $2 \times 2$ quantum minors appear – generate $H$-primes. Ten of the remaining cases – all except the one in position (2,2) – can be handled by similar means. Namely, in each of these cases, one can show that the factor of $A$ by the indicated ideal is an iterated skew polynomial algebra over the factor algebra $\mathcal{O}_q(M_2(k))/\langle\square\rangle$. Hence, these ten ideals are $H$-primes.

For the final case of Figure 3.1, we use the following result of Jordan.

**3.4. Proposition.** [12, Proposition 1] *Let $\sigma$ be an automorphism and $\delta$ a $\sigma$-derivation of a domain $S$. Let $T = S[x; \sigma, \delta]$. Let $c$ be a normal element of $T$ of the form $dx + e$, where $d, e \in S$ and $d \neq 0$. Let $\beta$ be the automorphism of $T$ such that $ct = \beta(t)c$ for all $t \in T$. Then $\beta(S) = S$, the element $d$ is normal in $S$, and $\beta(s)d = d\sigma(s)$ for all $s \in S$. Furthermore, if $e$ is regular modulo the ideal $Sd = dS$, then $T/Tc$ is a domain.* $\square$

**3.5.** To complete the cases indicated in Figure 3.1, we prove that the ideal $P$ of $A$ generated by the elements $u = [12 \mid 23]$ and $v = [23 \mid 12]$ is completely prime. Recall (e.g., from [**15**, Lemma 4.5.1]) that $u$ and $v$ scalar-commute with all the $X_{ij} \in A$ (that is, there are scalars $\lambda_{ij}$ and $\mu_{ij}$ such that $uX_{ij} = \lambda_{ij}X_{ij}u$ and $vX_{ij} = \mu_{ij}X_{ij}v$ for all $i,j$), whence $u$ and $v$ are normal elements of any subalgebra of $A$ which contains them and is generated by some of the $X_{ij}$.

Write $A = T[X_{33}; \tau_{33}, \delta_{33}]$, where $T = k\langle X_{11}, \ldots, X_{32}\rangle$, and set $Q = uT + vT$. Since $u$ and $v$ scalar-commute with $X_{33}$, they are $\delta_{33}$-constants, and so $Q$ is a $(\tau_{33}, \delta_{33})$-ideal of $T$. Consequently, $P = QA$ and $A/P \cong (T/Q)[X_{33}; \tau_{33}, \delta_{33}]$. Thus, it suffices to show that $T/Q$ is a domain.

Next, write $T = R[X_{32}; \tau_{32}, \delta_{32}]$, where $R = k\langle X_{11}, \ldots, X_{31}\rangle$. We can view $R$ as an iterated skew polynomial ring of the form $R_0[X_{31}; \sigma_{31}][X_{21}; \sigma_{21}][X_{11}; \sigma_{11}, \delta_{11}]$, where $R_0 = k\langle X_{12}, X_{13}, X_{22}, X_{23}\rangle$. Note that $R_0$ is a copy of $\mathcal{O}_q(M_2(k))$, in which $u$ corresponds to the quantum determinant, and so $R_0/uR_0$ is a domain. Hence, the algebra $S = R/uR$, which is an iterated skew polynomial extension of $R_0/uR_0$, is a domain. Since $u$ scalar-commutes with $X_{32}$, the ideal $uR$ is a $(\tau_{32}, \delta_{32})$-ideal of $R$, and we may identify the algebra $T' = T/uT$ with a skew polynomial ring $S[X_{32}; \tau_{32}, \delta_{32}]$. Set $c = v + uT$, so that $c$ is a normal element of $T'$ and $T/Q \cong T'/cT'$. Thus, we need to show that $T'/cT'$ is a domain.

Now $c = dX_{32} + e$ where $d = X_{21} + uT$ and $e = -qX_{22}X_{31} + uT$. It is clear that $d$ is a nonzero, normal element of $S$. If we show that $e$ is regular modulo $dS$, then Proposition 3.4 will imply that $T'/cT'$ is a domain. Since $X_{21}$ scalar-commutes with $X_{ij}$ for $(i,j) \leq_{\text{lex}} (3,1)$, we see that

$$S/dS \cong R/(uR + X_{21}R) \cong (R_0/uR_0)[X_{31}; \sigma_{31}][X_{11}; \sigma_{11}],$$

and so $S/dS$ is a domain. Further, $X_{22} \notin uR_0$, whence $e \notin dS$, and so $e$ is indeed regular modulo $dS$. This concludes the proof that $uA + vA$ is a completely prime ideal of $A$. It is obviously an $H$-ideal, and therefore an $H$-prime.

**3.6.** Steps 1 and 2 of (3.1) are now complete. We deal with Step 3 as follows.

Let $B$ denote the subalgebra $k\langle X_{12}, X_{13}, X_{22}, X_{23}\rangle \subseteq A$; thus $B$ is a copy of $\mathcal{O}_q(M_2(k))$. Observe that there is a $k$-algebra retraction $\theta : A \to B$ such that $\theta$ sends $X_{ij} \mapsto 0$ when $i = 3$ or $j = 1$ and fixes $X_{ij}$ otherwise. In particular, since $\theta$ is surjective, it maps ideals of $A$ onto ideals of $B$. Now observe that if $P$ is an $H$-prime of $A$ generated by one of the schematics in Figure 3.1, the ideal $\theta(P)$ of $B$ depends only on the column in which $P$ occurs. Specifically, column by column the $\theta(P)$ are given by

$$\begin{matrix} \circ\circ \\ \circ\circ \end{matrix} \quad \square \quad \begin{matrix} \circ\bullet \\ \circ\circ \end{matrix} \quad \begin{matrix} \bullet\bullet \\ \circ\circ \end{matrix} \quad \begin{matrix} \circ\bullet \\ \circ\bullet \end{matrix} \quad \begin{matrix} \bullet\bullet \\ \circ\bullet \end{matrix}$$

respectively. Since these 6 ideals of $B$ are distinct, we conclude that any two $H$-primes of $A$ arising from different columns of Figure 3.1 must be distinct.

A symmetric argument shows that $H$-primes arising from different rows of the figure are distinct, and therefore the $H$-primes generated by the 36 schematics in Figure 3.1 are all distinct. Since we have already seen that $A$ has exactly 36 $H$-primes not containing the quantum determinant, we conclude that Figure 3.1 completely describes this portion of $H$-$\operatorname{spec} A$.

## 4. Start of the final case

We continue to focus on $A = \mathcal{O}_q(M_3(k))$ in this section.

**4.1.** It remains to determine the $H$-primes of $A$ which contain the quantum determinant but do not contain all $2 \times 2$ quantum minors; these are the $H$-primes which lie in the sets $\operatorname{spec}_{\boldsymbol{r},\boldsymbol{c}} A$ for $\boldsymbol{r}, \boldsymbol{c} = (1,2), (1,3), (2,3)$. As indicated in Figure 2.3, $R_{12}^-$, $R_{13}^-$, and $R_{23}^-$ together have exactly 12 $H$-primes, and similarly for the corresponding $R_{\boldsymbol{r}}^+$. Hence, Theorem 1.6 shows that there are exactly 144 $H$-primes in $A$ containing the quantum determinant but not all $2 \times 2$ quantum minors. To determine them, we again follow the method of (3.1). This time, Step 1 involves all choices of $Q^+$ in $H$-$\operatorname{spec} R_{\boldsymbol{r}}^+$ and $Q^-$ in $H$-$\operatorname{spec} R_{\boldsymbol{c}}^-$ as $\boldsymbol{r}$ and $\boldsymbol{c}$ run through $(1,2), (1,3), (2,3)$. We obtain 144 generating sets which we display in Figure 4.1 below; here the diamond which appears in four cases stands for the quantum determinant of $A$.

**4.2.** Lemma 3.2 shows that 88 of the displayed schematics in Figure 4.1 – those in which no $2 \times 2$ or $3 \times 3$ quantum minors appear – generate $H$-primes. Among the remaining cases, there are 46 for which $A$ modulo the indicated ideal is isomorphic to an iterated skew polynomial extension of either $\mathcal{O}_q(M_2(k))/\langle\square\rangle$ or $\mathcal{O}_q(M_{2,3}(k))/\langle\square\square\rangle$. Hence, in these cases too the schematics generate $H$-primes (recall from [**5**, Corollary 2.6] that $\langle\square\square\rangle$ is a completely prime ideal of $\mathcal{O}_q(M_{2,3}(k))$). For instance, the factor algebras corresponding to the second and fourth positions of the first row are isomorphic to iterated skew polynomial extensions of $\mathcal{O}_q(M_{2,3}(k))/\langle\square\square\rangle$ and $\mathcal{O}_q(M_2(k))/\langle\square\rangle$ with respect to the variables $X_{31}$, $X_{21}$, $X_{11}$ (in that order). To handle cases like the ninth and tenth positions of the second row, note that $A/\langle X_{12}, X_{22}, X_{32}\rangle \cong \mathcal{O}_q(M_{3,2}(k))$.

Ten cases remain, which we display in Figure 4.2. The first of these cases is the ideal generated by the quantum determinant, an ideal which we already know is $H$-prime [**5**, Theorem 2.5]. We can handle the second, third, and fourth cases using Proposition 3.4, as indicated in (4.3) below. However, the final 6 cases require different methods, analogous to the preferred basis method used in [**5**]. We develop these methods in Sections 5 and 6, and conclude our work with $\mathcal{O}_q(M_3(k))$ in Section 7.

**4.3.** We show that the second schematic in Figure 4.2 generates an $H$-prime of $A$; the third and fourth schematics can be analyzed in the same way. Thus, we want to show that the ideal $P$ of $A$ generated by $X_{13}$ and the quantum determinant is completely prime.

By Lemma 3.2, the algebra $T = A/X_{13}A$ is an 8-fold iterated skew polynomial extension of $k$, and we write $T = S[X_{33}; \tau_{33}, \delta_{33}]$ where $S$ is the $k$-subalgebra of $T$ generated by (the cosets of) $X_{11}, X_{12}, X_{21}, X_{22}, X_{23}, X_{31}, X_{32}$. Let $c = [123 \mid 123] + X_{13}A$, a central element of $T$. Using the $q$-Laplace expansion of $[123 \mid 123]$ along the third column (e.g., [**15**, Corollary 4.4.4]), we can write $c = dX_{33} + e$ where $d = [12 \mid 12]$ and $e = -q[13 \mid 12]X_{23}$ (viewed as elements of $S$). It is clear that $d$ is a nonzero, normal element of $S$. Observe that $S/dS$ is an iterated skew polynomial extension of $\mathcal{O}_q(M_2(k))/\langle\square\rangle$, hence a domain, and that the image of $e$ in this domain is nonzero. Thus $e$ is regular modulo $dS$, and so Proposition 3.4 shows that $T/Tc$ is a domain. Therefore $A/P$ is a domain, as desired.

Figure 4.1

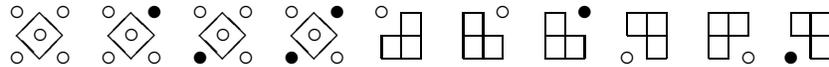

Figure 4.2

## 5. Special bases for $\mathcal{O}_q(M_n(k))$

We fix $A = \mathcal{O}_q(M_n(k))$ throughout this section, with $n$ any positive integer and $q$ an arbitrary nonzero scalar. A key ingredient in proving that quantum determinantal ideals in $A$ are completely prime [**5**, Theorem 2.5] was the *preferred basis* for $A$ developed in [**5**]. For our present calculations, we shall need that basis as well as three similar ones. We review the preferred basis and establish the others in this section.

**5.1.** We shall work with index sets contained in $\{1, \ldots, n\}$, and with the *row* and *column orderings* $\leq_r$ and $\leq_c$ on index sets defined in [**5**, (1.2)], as follows. Let $A, B \subseteq \{1, \ldots, n\}$, not necessarily of the same cardinality. First, write $A$ and $B$ in descending order:

$$A = \{a_1 > a_2 > \cdots > a_\alpha\} \qquad \text{and} \qquad B = \{b_1 > b_2 > \cdots > b_\beta\}.$$

Define $A \leq_r B$ to mean that $\alpha \geq \beta$ and $a_i \geq b_i$ for $i = 1, \ldots, \beta$. Now write $A$ and $B$ in ascending order:

$$A = \{a_1 < a_2 < \cdots < a_\alpha\} \qquad \text{and} \qquad B = \{b_1 < b_2 < \cdots < b_\beta\}.$$

Define $A \leq_c B$ to mean that $\alpha \geq \beta$ and $a_i \leq b_i$ for $i = 1, \ldots, \beta$.

By an *index pair* we will mean a pair $(I, J)$ where $I, J \subseteq \{1, \ldots, n\}$ and $|I| = |J|$. We will occasionally need to order index pairs by the relation $(\leq_r, \leq_c)$, that is, we define $(I, J) \leq (I', J')$ if and only if $I \leq_r I'$ and $J \leq_c J'$.

**5.2.** Our bases will be indexed by certain bitableaux (pairs of tableaux with a common shape). Recall that, in general, a *tableau* consists of a Young diagram with entries in each box. For present purposes, define an *allowable tableau* to be any tableau with entries from $\{1, \ldots, n\}$ such that there are no repeated entries in any row. Then define an *allowable bitableau* to be any bitableau $(T, T')$ where $T$ and $T'$ are allowable. This definition of allowable bitableaux is more general than that in [**5**, (1.3)], and so we must be careful when citing results from that paper.

Let us say that a tableau is *increasing* (respectively, *decreasing*) provided all its rows are strictly increasing (respectively, strictly decreasing) and all its columns are nondecreasing (respectively, nonincreasing).

Four kinds of bitableaux will be important. We say that a bitableau $(T, T')$ is

- *preferred* if $T$ is decreasing and $T'$ is increasing;
- *antipreferred* if $T$ is increasing and $T'$ is decreasing;
- *standard* if $T$ and $T'$ are increasing;
- *antistandard* if $T$ and $T'$ are decreasing.

All these types of bitableaux are allowable; those we have labelled 'preferred' are both allowable and preferred in the sense of [**5**]. Note that $(T, T')$ is antipreferred if and only if $(T', T)$ is preferred.

**5.3.** As in [**5**, (1.3)], we may indicate an allowable bitableau $(T, T')$ in the form

$$\begin{pmatrix} I_1 & J_1 \\ I_2 & J_2 \\ \vdots & \vdots \\ I_l & J_l \end{pmatrix}$$

where $(I_1, J_1), \ldots, (I_l, J_l)$ are index pairs such that $|I_1| \geq |I_2| \geq \cdots \geq |I_l|$. However, these index pairs do not determine $(T, T')$ unless we specify how the entries in the rows of $T$ and $T'$ are to be listed. The four special types of allowable bitableaux defined in (5.2) correspond to the following conditions on the index sets $I_j$ and $J_j$:

- preferred: $I_1 \leq_r I_2 \leq_r \cdots \leq_r I_l$ and $J_1 \leq_c J_2 \leq_c \cdots \leq_c J_l$;
- antipreferred: $I_1 \leq_c I_2 \leq_c \cdots \leq_c I_l$ and $J_1 \leq_r J_2 \leq_r \cdots \leq_r J_l$;
- standard: $I_1 \leq_c I_2 \leq_c \cdots \leq_c I_l$ and $J_1 \leq_c J_2 \leq_c \cdots \leq_c J_l$;
- antistandard: $I_1 \leq_r I_2 \leq_r \cdots \leq_r I_l$ and $J_1 \leq_r J_2 \leq_r \cdots \leq_r J_l$.

As in [**5**, (1.4)], we define

$$[T \mid T'] = [I_1 \mid J_1][I_2 \mid J_2] \cdots [I_l \mid J_l].$$

The element $[T \mid T']$ is called a *preferred* (respectively, *antipreferred*, *standard*, *antistandard*) *product* if the bitableau $(T, T')$ is preferred (respectively, antipreferred, standard, antistandard). These four types of products of quantum minors will make up four bases for $A$. In particular:

**Theorem.** [**5**, Corollary 1.10] *The preferred products $[T \mid T']$, as $(T, T')$ runs over all preferred bitableaux, form a basis for $A$.* □

**5.4.** Let $\boldsymbol{\tau}$ denote the "transpose" automorphism of $A$, that is, the $k$-algebra automorphism such that $\boldsymbol{\tau}(X_{ij}) = X_{ji}$ for all $i, j$. By [**15**, Lemma 4.3.1], $\boldsymbol{\tau}([I \mid J]) = [J \mid I]$ for all index pairs $(I, J)$, and consequently

$$\boldsymbol{\tau}([T \mid T']) = [T' \mid T]$$

for all allowable bitableaux $(T, T')$. In particular, $\boldsymbol{\tau}$ maps the set of preferred products bijectively onto the set of antipreferred products. Thus, the previous theorem immediately implies the following result:

**Theorem.** *The antipreferred products $[T \mid T']$, as $(T, T')$ runs over all antipreferred bitableaux, form a basis for $A$.* □

**5.5.** To facilitate working with products of quantum minors corresponding to the various types of bitableaux we have defined, we label several aspects of tableaux and bitableaux as follows.

For any allowable tableau $T$, let $\mu(T)$ and $\mu'(T)$ denote the tableaux with the same shape as $T$ and entries as follows: each row of length $l$ is filled $1, 2, \ldots, l$ in $\mu(T)$ and is filled $l, l-1, \ldots, 1$ in $\mu'(T)$. Thus $\mu(T)$ is increasing, and $\mu'(T)$ is decreasing. Note that $\mu(T)$ and $\mu'(T)$ are based on the the same list of index sets; in particular, $[T \mid \mu(T)] = [T \mid \mu'(T)]$ and $[\mu(T) \mid T'] = [\mu'(T) \mid T']$ for any allowable bitableau $(T, T')$. For $l = 1, \ldots, n$, let $\rho_l(T)$ be the number of rows of $T$ of length $l$ or greater, and set $\overline{\rho}(T) = (\rho_1(T), \rho_2(T), \ldots, \rho_n(T))$.

Recall that there is a natural $\mathbb{Z}^n \times \mathbb{Z}^n$ bigrading on $A$, under which each $X_{ij}$ has bidegree $(\epsilon_i, \epsilon_j)$ where $\epsilon_1, \ldots, \epsilon_n$ is the standard basis for $\mathbb{Z}^n$. Any quantum minor $[I \mid J]$ is homogeneous of bidegree $(\sum_{i \in I} \epsilon_i, \sum_{j \in J} \epsilon_j)$. For any homogeneous element $x \in A$, label the bidegree of $x$ as

$$(\overline{r}(x), \overline{c}(x)) = (r_1(x), r_2(x), \ldots, r_n(x), c_1(x), c_2(x), \ldots, c_n(x)).$$

Thus, with respect to the usual PBW basis of ordered monomials, $r_l(x)$ records the number of $X_{l?}$ factors in each monomial in $x$, and $c_l(x)$ the number of $X_{?l}$ factors. If $(T, T')$ is an allowable bitableau, then $r_l[T \mid T']$ is the number of $l$'s in $T$ and $c_l[T \mid T']$ is the number of $l$'s in $T'$. Note that $c_l[T \mid \mu(T)] = r_l[\mu'(T) \mid T'] = \rho_l(T)$ for all $l$; in particular, $\overline{c}[T \mid \mu(T)] = \overline{\rho}(T)$.

We shall also need the notion of the *bicontent* of an allowable bitableau $(T, T')$, which is the pair of multisets

$$(1^{r_1} 2^{r_2} \cdots n^{r_n}, 1^{c_1} 2^{c_2} \cdots n^{c_n}),$$

where $r_i = r_i([T \mid T'])$ and $c_j = c_j([T \mid T'])$. In other words, the bicontent $(R, C)$ of $(T, T')$ is obtained by recording the entries of $T$ in the multiset $R$ and those of $T'$ in the multiset $C$.

We shall write $<_{\text{rlex}}$ for the reverse lexicographic order on $n$-tuples of integers.

**5.6. Lemma.** *If $(T, T')$ is an allowable bitableau, then*

$$\Delta[T \mid T'] = [T \mid \mu(T)] \otimes [\mu(T) \mid T'] + \sum_i X_i \otimes Y_i$$

*where the $X_i$ and $Y_i$ are homogeneous elements of $A$ with $\overline{c}(X_i) = \overline{r}(Y_i) >_{\text{rlex}} \overline{\rho}(T)$. More precisely, for each $i$ the terms $X_i$ and $Y_i$ are products of the forms $[T \mid S_i]$ and $[S_i \mid T']$, respectively.*

*Proof.* Let $(S, S')$ be the allowable bitableau such that $S$ has the same shape as $T$ and each row of $S$ (respectively, $S'$) contains the same entries as the corresponding row of $T$ (respectively, $T'$) but listed in descending (respectively, ascending) order. Then $(S, S')$ is allowable in the sense of [**5**, (1.3)]. Moreover, $[T \mid T'] = [S \mid S']$, and

$$[T \mid R] \otimes [R \mid T'] = [S \mid R] \otimes [R \mid S']$$

for any allowable tableau $R$ with the same shape as $T$, while also $\mu(T) = \mu(S)$ and $\overline{\rho}(T) = \overline{\rho}(S)$. Therefore the lemma follows from the case $t = n$ of [**5**, Lemma 2.3]. □

**5.7. Corollary.** *Let $x = \sum_{i=1}^{m} \alpha_i [T_i \mid T_i']$ where the $\alpha_i \in k$ and the $(T_i, T_i')$ are allowable bitableaux. Let $\overline{\rho}_{\min}$ be the minimum of the $\overline{\rho}(T_i)$ under reverse lexicographic order, and assume that there exists $m'$ such that $\overline{\rho}(T_i) = \overline{\rho}_{\min}$ for $i \leq m'$ and $\overline{\rho}(T_i) >_{\text{rlex}} \overline{\rho}_{\min}$ for $i > m'$. If $U$ is a homogeneous subspace of $A \otimes A$ and $\Delta(x) \in U$, then*

$$\sum_{i=1}^{m'} \alpha_i [T_i \mid \mu(T_i)] \otimes [\mu(T_i) \mid T_i'] \in U.$$

*Proof.* Applying Lemma 5.6 to each $\Delta[T_i \mid T_i']$ and collecting terms, we see that

$$\Delta(x) = \sum_{i=1}^{m'} \alpha_i [T_i \mid \mu(T_i)] \otimes [\mu(T_i) \mid T_i'] + \sum_j X_j \otimes Y_j$$

where the $X_j$ and $Y_j$ are homogeneous with $\overline{c}(X_j) >_{\text{rlex}} \overline{\rho}_{\min}$. Since $\overline{c}[T_i \mid \mu(T_i)] = \overline{\rho}_{\min}$ for $i \leq m'$, all of the $X_j$ belong to different homogeneous components of $A$ than the $[T_i \mid \mu(T_i)]$ for $i \leq m'$. Given that $U$ is homogeneous, the desired conclusion follows. $\square$

**5.8. Theorem.** *The standard products $[T \mid T']$, as $(T, T')$ runs over all standard bitableaux, form a basis for $A$.*

**Remark.** The case of this theorem when $q$ is not a root of unity follows from work of Huang and Zhang [**11**, Theorem 9].

*Proof.* We begin with linear independence. Let $x = \sum_{i=1}^{m} \alpha_i [T_i \mid T_i']$ where the $\alpha_i$ are nonzero scalars and the $(T_i, T_i')$ are distinct standard bitableaux. To show that $x \neq 0$, we will show that $\Delta(x) \neq 0$.

Let $\overline{\rho}_{\min}$ be the minimum of the $\overline{\rho}(T_i)$ under reverse lexicographic order. Without loss of generality, we may assume that there exists $m'$ such that $\overline{\rho}(T_i) = \overline{\rho}_{\min}$ for $i \leq m'$ and $\overline{\rho}(T_i) >_{\text{rlex}} \overline{\rho}_{\min}$ for $i > m'$. If $\Delta(x) = 0$, then

$$\sum_{i=1}^{m'} \alpha_i [T_i \mid \mu'(T_i)] \otimes [\mu'(T_i) \mid T_i'] = 0$$

by Corollary 5.7. (Recall that $[T_i \mid \mu(T_i)] = [T_i \mid \mu'(T_i)]$ and $[\mu(T_i) \mid T_i'] = [\mu'(T_i) \mid T_i']$.) For $1 \leq i < j \leq m$, either $T_i \neq T_j$ or $T_i' \neq T_j'$, whence either $(T_i, \mu'(T_i)) \neq (T_j, \mu'(T_j))$ or $(\mu'(T_i), T_i') \neq (\mu'(T_j), T_j')$. Moreover, the $(T_i, \mu'(T_i))$ are antipreferred bitableaux and the $(\mu'(T_i), T_i')$ are preferred. Since the antipreferred (respectively, preferred) products in $A$ are linearly independent, the tensors $[T_i \mid \mu'(T_i)] \otimes [\mu'(T_i) \mid T_i']$ in $A \otimes A$ must be linearly independent. But then $\alpha_i = 0$ for $1 \leq i \leq m'$, which contradicts our assumptions. Thus $\Delta(x) \neq 0$, whence $x \neq 0$.

Therefore the standard products in $A$ are linearly independent.

To see that the standard products span $A$, it suffices to show that they span each homogeneous component. Let $V$ be the homogeneous component of $A$ of a given bidegree $\delta$, and let $\gamma$ denote the corresponding bicontent.

Write $A$ as an iterated skew polynomial extension of $k$ as in (1.1), with the generators $X_{ij}$ in lexicographic order. The corresponding PBW basis for $A$ consists of monomials $X_{i_1 j_1} X_{i_2 j_2} \cdots X_{i_r j_r}$ such that $i_1 \leq i_2 \leq \cdots \leq i_r$ and $j_l \leq j_{l+1}$ whenever $i_l = j_l$. Those monomials of this type with bidegree $\delta$ form a basis for $V$. Hence, $\dim_k V$ equals the number of two-rowed matrices

$$\begin{pmatrix} i_1 & i_2 & \cdots & i_r \\ j_1 & j_2 & \cdots & j_r \end{pmatrix}$$

with

- entries from $\{1, \ldots, n\}$;
- columns in lexicographic order;
- $(\{i_1, \ldots, i_r\}, \{j_1, \ldots, j_r\}) = \gamma$ (as pairs of multisets).

By the Robinson-Schensted-Knuth Theorem [**3**, p. 40], the set of these matrices is in bijection with the set of those bitableaux $(Q, P)$ of bicontent $\gamma$ having nondecreasing rows and strictly increasing columns. The latter conditions just say that $(Q^{\mathrm{tr}}, P^{\mathrm{tr}})$ is a standard bitableau in our sense. Hence, $\dim_k V$ equals the number of standard bitableaux $(T, T')$ with bicontent $\gamma$. Since the corresponding standard products $[T \mid T']$ are linearly independent elements of $V$, they must span $V$.

Therefore $A$ is spanned by standard products. $\square$

**5.9. Corollary.** *The antistandard products $[T \mid T']$, as $(T, T')$ runs over all antistandard bitableaux, form a basis for $A$.*

*Proof.* Let $\boldsymbol{\rho} = \rho_q$ denote the $k$-algebra anti-automorphism of $A$ such that $\boldsymbol{\rho}(X_{ij}) = X_{n+1-j,\, n+1-i}$ for all $i, j$ [**15**, Proposition 3.7.1(3)]. Set $A' = \mathcal{O}_{q^{-1}}(M_n(k))$ with generators $X'_{ij}$, and write $[I \mid J]'$ for quantum minors in $A'$. There is a $k$-algebra anti-isomorphism $\boldsymbol{\alpha} : A \to A'$ such that $\boldsymbol{\alpha}(X_{ij}) = X'_{ji}$ for all $i, j$ [**15**, Remark 3.7.2]. Hence, $\boldsymbol{\delta} = \boldsymbol{\alpha}\boldsymbol{\rho} : A \to A'$ is a $k$-algebra isomorphism such that $\boldsymbol{\delta}(X_{ij}) = X'_{n+1-i,n+1-j}$ for all $i, j$. We claim that $\boldsymbol{\delta}([I \mid J]) = [\omega_0 I \mid \omega_0 J]'$ for all index pairs $(I, J)$, where $\omega_0 = \begin{pmatrix} 1 & 2 & \cdots & n \\ n & n-1 & \cdots & 1 \end{pmatrix}$ denotes the longest element of the symmetric group $S_n$.

Given an index pair $(I, J)$, consider the subalgebras
$$A_{I,J} = k\langle X_{ij} \mid i \in I,\, j \in J \rangle \subseteq A$$
$$A'_{\omega_0 I, \omega_0 J} = k\langle X'_{ij} \mid i \in \omega_0 I,\, j \in \omega_0 J \rangle \subseteq A',$$
and observe that $\boldsymbol{\delta}$ maps $A_{I,J}$ isomorphically onto $A'_{\omega_0 I, \omega_0 J}$. Write $I$ and $J$ in ascending order, say $I = \{i_1 < \cdots < i_t\}$ and $J = \{j_1 < \cdots < j_t\}$; then
$$\omega_0 I = \{n + 1 - i_t < n + 1 - i_{t-1} < \cdots < n + 1 - i_1\},$$
and similarly for $\omega_0 J$. There are $k$-algebra isomorphisms
$$\phi : \mathcal{O}_q(M_t(k)) \longrightarrow A_{I,J} \qquad \text{and} \qquad \phi' : \mathcal{O}_{q^{-1}}(M_t(k)) \longrightarrow A'_{\omega_0 I, \omega_0 J}$$
such that $\phi(X_{rc}) = X_{i_r j_c}$ and $\phi'(X'_{rc}) = X'_{n+1-i_{t+1-r},\, n+1-j_{t+1-c}}$ for $r, c = 1, \ldots, t$. The restriction of $\boldsymbol{\delta}$ to $A_{I,J}$ can now be factored as the composition
$$A_{I,J} \xrightarrow{\phi^{-1}} \mathcal{O}_q(M_t(k)) \xrightarrow{\boldsymbol{\delta}_t} \mathcal{O}_{q^{-1}}(M_t(k)) \xrightarrow{\phi'} A'_{\omega_0 I, \omega_0 J}$$
where $\boldsymbol{\delta}_t$ is the $k$-algebra isomorphism sending $X_{rc} \mapsto X'_{t+1-r,t+1-c}$ for $r, c = 1, \ldots, t$. By [**15**, Remark 4.1.2], $\boldsymbol{\delta}_t$ sends the quantum determinant of $\mathcal{O}_q(M_t(k))$ to the quantum determinant of $\mathcal{O}_{q^{-1}}(M_t(k))$, from which we conclude that $\boldsymbol{\delta}$ sends $[I \mid J]$ to $[\omega_0 I \mid \omega_0 J]'$, as claimed.

In view of the claim, $\boldsymbol{\delta}([T \mid T']) = [\omega_0 T \mid \omega_0 T']'$ for all allowable bitableaux $(T, T')$. It follows that $\boldsymbol{\delta}$ maps the set of antistandard products in $A$ bijectively onto the set of standard products in $A'$. Therefore the corollary follows from Theorem 5.8. $\square$

## 6. New $H$-primes in $\mathcal{O}_q(M_n(k))$

With the help of the bases constructed in the previous section, we identify some new $H$-primes in the algebra $A = \mathcal{O}_q(M_n(k))$, generated by quantum minors of two different sizes; here $n$ is again an arbitrary positive integer. For most of the section, $q$ can be an arbitrary nonzero scalar, but our proof of the main result requires $q$ to be a non-root of unity. As with the quantum determinantal ideals investigated in [**5**], we proceed by identifying certain inverse images under the comultiplication $\Delta : A \to A \otimes A$. First, we upgrade Corollary 5.7 to allow certain permutations of the tableaux $\mu(T)$, as follows.

**6.1.** Let $T$ be an allowable tableau, and $\sigma \in S_n$. Let $\mu_\sigma(T)$ denote the increasing tableau with the same shape as $T$ and entries as follows: each row of length $l$ in $\mu_\sigma(T)$ is filled with the numbers $\sigma(1), \sigma(2), \ldots, \sigma(l)$ in ascending order. In particular, $\mu_{\mathrm{id}}(T) = \mu(T)$.

**6.2. Lemma.** *Let $\sigma \in S_n$, and let $x = \sum_{i=1}^m \alpha_i [T_i \mid T_i']$ where the $\alpha_i \in k$ and the $(T_i, T_i')$ are allowable bitableaux. Let $\overline{\rho}_{\min}$ be the minimum of the $\overline{\rho}(T_i)$ under reverse lexicographic order, and assume that there exists $m'$ such that $\overline{\rho}(T_i) = \overline{\rho}_{\min}$ for $i \leq m'$ and $\overline{\rho}(T_i) >_{\mathrm{rlex}} \overline{\rho}_{\min}$ for $i > m'$. If $U$ is a homogeneous subspace of $A \otimes A$ and $\Delta(x) \in U$, then*

$$\sum_{i=1}^{m'} \alpha_i [T_i \mid \mu_\sigma(T_i)] \otimes [\mu_\sigma(T_i) \mid T_i'] \;\in\; U.$$

*Proof.* We first need to write each $\Delta[T_i \mid T_i']$ in a standard format. We have

$$[T_i \mid T_i'] = [I_{i1} \mid J_{i1}][I_{i2} \mid J_{i2}] \cdots [I_{is_i} \mid J_{is_i}]$$

where the $(I_{il}, J_{il})$ are index pairs. Since

$$\Delta[I_{il} \mid J_{il}] = \sum_{\substack{K \subseteq \{1,\ldots,n\} \\ |K|=|I_{il}|}} [I_{il} \mid K] \otimes [K \mid J_{il}]$$

for all $i, l$, we see that

$$(*) \qquad \Delta[T_i \mid T_i'] = \sum_j [T_i \mid S_{ij}] \otimes [S_{ij} \mid T_i']$$

where the $S_{ij}$ are allowable tableaux such that

$$(T_i, S_{ij}) = \begin{pmatrix} I_{i1} & K_{ij1} \\ I_{i2} & K_{ij2} \\ \vdots & \vdots \\ I_{is_i} & K_{ijs_i} \end{pmatrix} \qquad \text{and} \qquad (S_{ij}, T_i') = \begin{pmatrix} K_{ij1} & J_{i1} \\ K_{ij2} & J_{i2} \\ \vdots & \vdots \\ K_{ijs_i} & J_{is_i} \end{pmatrix}$$

for some index sets $K_{ijl}$ with $|K_{ijl}| = |I_{il}|$. In order to ensure that $S_{ij}$ is uniquely determined by its sequence of index sets, let us require $S_{ij}$ to be increasing. Thus, we can say

that in $(*)$, the $S_{ij}$ run bijectively through the increasing tableaux with the same shape as $T_i$.

Now because of $(*)$, we have

$$(**) \qquad \Delta(x) = \sum_{i=1}^{m} \sum_{j} \alpha_i [T_i \mid S_{ij}] \otimes [S_{ij} \mid T_i'].$$

For each $i$, there is an index, which we may take to be $j = 0$, such that $S_{i0} = \mu(T_i)$. By Lemma 5.6, $\overline{c}[T_i \mid S_{ij}] >_{\text{rlex}} \overline{\rho}(T_i)$ for all $i$ and all $j \neq 0$. As noted in (5.5), we also have $\overline{c}[T_i \mid S_{i0}] = \overline{\rho}(T_i)$. Hence,

$$(\dagger) \qquad \overline{c}[T_i \mid S_{ij}] >_{\text{rlex}} \overline{\rho}_{\min} \qquad \text{whenever } i > m' \text{ or } S_{ij} \neq \mu(T_i).$$

Let $(r_1, r_2, \ldots, r_n) = \overline{\rho}_{\min}$, and set $v = (r_{\tau(1)}, r_{\tau(2)}, \ldots, r_{\tau(n)})$ where $\tau = \sigma^{-1}$. We claim that $\overline{c}[T_i \mid S_{ij}] = v$ if and only if $i \leq m'$ and $S_{ij} = \mu_\sigma(T_i)$.

First let $i \leq m'$, so that $\overline{\rho}(T_i) = \overline{\rho}_{\min}$. In $\mu_\sigma(T_i)$, each $\sigma(l)$ appears precisely in those rows of length at least $l$, and so $c_{\sigma(l)}[T_i \mid \mu_\sigma(T_i)] = \rho_l(T_i) = r_l$. Hence, $c_l[T_i \mid \mu_\sigma(T_i)] = r_{\tau(l)}$ for all $l$, and thus $\overline{c}[T_i \mid \mu_\sigma(T_i)] = v$.

Conversely, let $i$ and $j$ be indices such that $\overline{c}[T_i \mid S_{ij}] = v$. Let $\widetilde{S}_{ij}$ be the increasing tableau with the same shape as $S_{ij}$ in which the entries of any given row are obtained by applying $\tau$ to the entries of the corresponding row of $S_{ij}$. For $1 \leq l \leq n$, observe that $\tau(l)$ appears in $\widetilde{S}_{ij}$ exactly as many times as $l$ appears in $S_{ij}$, whence $c_{\tau(l)}[T_i \mid \widetilde{S}_{ij}] = c_l[T_i \mid S_{ij}] = r_{\tau(l)}$. Thus $c_l[T_i \mid \widetilde{S}_{ij}] = r_l$ for all $l$, that is, $\overline{c}[T_i \mid \widetilde{S}_{ij}] = \overline{\rho}_{\min}$. In view of $(\dagger)$, we must have $i \leq m'$ and $\widetilde{S}_{ij} = \mu(T_i)$. Now all rows of length $l$ in $\widetilde{S}_{ij}$ have entries $1, 2, \ldots, l$. Hence, the entries in any row of length $l$ in $S_{ij}$ are $\sigma(1), \sigma(2), \ldots, \sigma(l)$, and therefore $S_{ij} = \mu_\sigma(T_i)$. This completes the proof of the claim.

Because of the claim, the terms $\alpha_i [T_i \mid \mu_\sigma(T_i)] \otimes [\mu_\sigma(T_i) \mid T_i']$ with $i \leq m'$ in $(**)$ lie in different homogeneous components of $A \otimes A$ than all the other terms. Since $U$ is homogeneous, the lemma follows. $\square$

**6.3. Proposition.** Let $V_1, W_1, \ldots, V_p, W_p$ be $k$-subspaces of $A$, and let

$$x = \sum_{i=1}^{m} \alpha_i [T_i \mid T_i'] \;\; \in \;\; \Delta^{-1}\Big(\sum_{l=1}^{p} V_l \otimes W_l\Big),$$

where the $\alpha_i$ are nonzero scalars and the $(T_i, T_i')$ are distinct allowable bitableaux. Let $\overline{\rho}_{\min}$ be the minimum of the $n$-tuples $\overline{\rho}(T_i)$ under reverse lexicographic order. Assume one of the following:

(1) Each $(T_i, T_i')$ is standard, either each $V_l$ is spanned by standard products or each $V_l$ is spanned by antipreferred products, and either each $W_l$ is spanned by standard products or each $W_l$ is spanned by preferred products.
(2) Each $(T_i, T_i')$ is antistandard, either each $V_l$ is spanned by antistandard products or each $V_l$ is spanned by preferred products, and either each $W_l$ is spanned by antistandard products or each $W_l$ is spanned by antipreferred products.

(3) Each $(T_i, T_i')$ is preferred, either each $V_l$ is spanned by preferred products or each $V_l$ is spanned by antistandard products, and either each $W_l$ is spanned by preferred products or each $W_l$ is spanned by standard products.

(4) Each $(T_i, T_i')$ is antipreferred, either each $V_l$ is spanned by antipreferred products or each $V_l$ is spanned by standard products, and either each $W_l$ is spanned by antipreferred products or each $W_l$ is spanned by antistandard products.

Then for each index $i$ such that $\overline{\rho}(T_i) = \overline{\rho}_{\min}$ and each $\sigma \in S_n$, there is an index $l(i, \sigma)$ such that $[T_i \mid \mu_\sigma(T_i)] \in V_{l(i,\sigma)}$ and $[\mu_\sigma(T_i) \mid T_i'] \in W_{l(i,\sigma)}$.

*Proof.* Fix $\sigma \in S_n$. Note that the $V_l$ and $W_l$ are homogeneous subspaces of $A$. Without loss of generality, there exists $m'$ such that $\overline{\rho}(T_i) = \overline{\rho}_{\min}$ for $i \leq m'$ and $\overline{\rho}(T_i) >_{\text{rlex}} \overline{\rho}_{\min}$ for $i > m'$. By Lemma 6.2,

$$(*) \qquad \sum_{i=1}^{m'} \alpha_i [T_i \mid \mu_\sigma(T_i)] \otimes [\mu_\sigma(T_i) \mid T_i'] \in \sum_l V_l \otimes W_l.$$

Note that for $i_1 \neq i_2$, we must have either $(T_{i_1}, \mu_\sigma(T_{i_1})) \neq (T_{i_2}, \mu_\sigma(T_{i_2}))$ or $(\mu_\sigma(T_{i_1}), T_{i_1}') \neq (\mu_\sigma(T_{i_2}), T_{i_2}')$.

In case (1), we can rearrange the entries in the rows of the $\mu_\sigma(T_i)$ to be increasing or decreasing as desired to match the given subcases of (1). Hence, we can arrange things so that the $(T_i, \mu_\sigma(T_i))$ are either all standard or all antipreferred, and so that the $(\mu_\sigma(T_i), T_i')$ are either all standard or all preferred. Correspondingly, there is a basis for $A \otimes A$ consisting of all tensors $[R \mid R'] \otimes [S \mid S']$ where the $(R, R')$ are either all standard or all antipreferred, and the $(S, S')$ are either all standard or all preferred. Then $\sum_l V_l \otimes W_l$ is spanned by a subset of this basis. Thus, it follows from $(*)$ that for $1 \leq i \leq m'$, there is some $l(i)$ such that $[T_i \mid \mu_\sigma(T_i)] \in V_{l(i)}$ and $[\mu_\sigma(T_i) \mid T_i'] \in W_{l(i)}$.

Cases (2), (3), and (4) are handled analogously. □

Although for the present section we shall only make use of Proposition 6.3 in the case where $\sigma = \text{id}$ (see the proof of Theorem 6.6), the general case will be required in (7.2).

**6.4.** Recall from (5.1) the ordering $(\leq_r, \leq_c)$, written $\leq$, on index pairs. An index pair $(I, J)$ with $|I| = |J| = t$ will be called a $t \times t$ *index pair*. A set $\mathcal{P}$ of $t \times t$ index pairs is *hereditary* among such pairs provided that for any $t \times t$ index pairs $(I', J') \leq (I, J)$, if $(I, J) \in \mathcal{P}$ then $(I', J') \in \mathcal{P}$.

**Lemma.** *Let $1 \leq t \leq n$, let $\mathcal{P}$ be a set of $t \times t$ index pairs which is hereditary among $t \times t$ index pairs, and set*

$$P = \langle [I \mid J] \mid |I| = t+1 \text{ or } (I, J) \in \mathcal{P} \rangle.$$

*Then $P$ has a basis consisting of those preferred products $[T \mid T']$ such that either $T$ has more than $t$ columns or $(T, T')$ has a row that lies in $\mathcal{P}$.*

*Proof.* It is clear that $P$ contains all preferred products of the described types, so we just need to show that these preferred products span $P$. We first consider the case when $\mathcal{P}$ is empty. What is needed is to show that the ideal $I_t := \langle [I \mid J] \mid |I| = t+1 \rangle$ is spanned by those preferred products $[T \mid T']$ such that $T$ has more than $t$ columns. We do so by a downward induction on $t$. If $t = n$, then $I_t = 0$ and the result holds trivially.

Now suppose that $t < n$ and that $I_{t+1}$ is spanned by the preferred products $[T \mid T']$ such that $T$ has at least $t+2$ columns. In view of [5, Corollary A.2], the set $\sum_{|I|=t+1}[I \mid J]A$ is a two-sided ideal of $A$, whence

$$I_t = \sum_{|I|=t+1} [I \mid J]A.$$

Consider one of the right ideals $[I \mid J]A$ in the sum above, where $(I, J)$ is a $(t+1) \times (t+1)$ index pair. Writing elements from $A$ using the preferred basis, we see that $[I \mid J]A$ is spanned by terms $[I \mid J][R \mid R']$ where $[R \mid R']$ is a preferred product. If $R$ has more than $t+1$ columns, then $[I \mid J][R \mid R']$ lies in $I_{t+1}$, and so it is a linear combination of preferred products of the desired type. Otherwise, $[I \mid J][R \mid R'] = [S \mid S']$ where $(S, S')$ is an allowable bitableau and $S$ has $t+1$ columns. By [5, Corollary 1.8], $[S \mid S']$ is a linear combination of preferred products $[T_i \mid T'_i]$ where each $T_i$ has at least $t+1$ columns. This shows that $I_t$ is spanned by those preferred products $[T \mid T']$ such that $T$ has at least $t+1$ columns, and completes the induction. Thus, we have shown that the lemma holds when $\mathcal{P}$ is empty.

In the general case, with another application of [5, Corollary A.2] we find that

$$P = I_t + \sum_{(I,J) \in \mathcal{P}} [I \mid J]A.$$

Let $(I, J) \in \mathcal{P}$, and let $[R \mid R']$ be an arbitrary preferred product in $A$. If $R$ has more than $t$ columns, then $[I \mid J][R \mid R']$ lies in $I_t$, and so it is a linear combination of preferred products $[T \mid T']$ such that $T$ has at least $t+1$ columns. Otherwise, $[I \mid J][R \mid R'] = [S \mid S']$ where $(S, S')$ is an allowable bitableau and $S$ has $t$ columns. By [5, Corollary 1.8], $[S \mid S']$ is a linear combination of preferred products $[T_i \mid T'_i]$ with top rows $(X_i, Y_i)$ such that either $|X_i| > t$ or $(X_i, Y_i) \leq (I, J)$. In the first case, $T_i$ has at least $t+1$ columns, while in the second, $(X_i, Y_i) \in \mathcal{P}$. Therefore $P$ is spanned by preferred products of the desired types. $\square$

**6.5. Lemma.** *Let $1 \leq t \leq n$. The ideal $\langle X_{ij} \mid j > t \rangle$ (respectively, $\langle X_{ij} \mid i > t \rangle$) has a basis consisting of those preferred products containing a row $(I, J)$ with $J \not\subseteq \{1, \ldots, t\}$ (respectively, $I \not\subseteq \{1, \ldots, t\}$). These ideals also have bases consisting of those antipreferred products containing rows of the given types.*

*Proof.* Let $U$ be the subspace of $A$ spanned by those preferred products containing a row $(I, J)$ with $J \not\subseteq \{1, \ldots, t\}$, and note that the remaining preferred products in $A$ span the subalgebra $\mathcal{O}_q(M_{n,t}(k))$. Thus $U \subseteq \langle X_{ij} \mid j > t \rangle$ and $A = U \oplus \mathcal{O}_q(M_{n,t}(k))$. Since the natural retraction of $A$ onto $\mathcal{O}_q(M_{n,t}(k))$ has kernel $\langle X_{ij} \mid j > t \rangle$, it follows that $U = \langle X_{ij} \mid j > t \rangle$, as desired.

The other three cases are proved in the same fashion. $\square$

The main theorem of this section, which we now prove, identifies a new supply of $H$-primes in $\mathcal{O}_q(M_n(k))$.

**6.6. Theorem.** *Assume that $q$ is not a root of unity. Let $1 \leq t, l, l' \leq n$, and set*

$$P_{t,l}^- = \langle [I \mid J] \mid |I| = t+1 \rangle + \langle [I \mid J] \mid |I| = t \text{ and } I \subseteq \{1, \ldots, l\} \rangle$$
$$P_{t,l}^+ = \langle [I \mid J] \mid |I| = t+1 \rangle + \langle [I \mid J] \mid |I| = t \text{ and } I \subseteq \{n-l+1, \ldots, n\} \rangle$$
$$Q_{t,l'}^- = \langle [I \mid J] \mid |I| = t+1 \rangle + \langle [I \mid J] \mid |I| = t \text{ and } J \subseteq \{1, \ldots, l'\} \rangle$$
$$Q_{t,l'}^+ = \langle [I \mid J] \mid |I| = t+1 \rangle + \langle [I \mid J] \mid |I| = t \text{ and } J \subseteq \{n-l'+1, \ldots, n\} \rangle.$$

*Then the ideals $P_{t,l}^\pm + Q_{t,l'}^\pm$ are $H$-prime ideals of $A$.*

**Remark.** The case $P_{t,l}^+ + Q_{t,l}^-$ of the theorem has been obtained by Lenagan and Rigal [**13**, Theorem 2] using different methods.

*Proof.* Set $V_1^\pm = P_{t,l}^\pm + \langle X_{ij} \mid j > t \rangle$ and $W_2^\pm = Q_{t,l'}^\pm + \langle X_{ij} \mid i > t \rangle$. We observe that these are $H$-prime ideals of $A$; for instance, $A/V_1^-$ is isomorphic to an iterated skew polynomial extension of the domain

$$\mathcal{O}_q(M_{l,t}(k))/\langle \text{all } t \times t \text{ quantum minors} \rangle.$$

By [**6**, (3.1) and Lemma 3.2], the ideals $(V_1^\pm \otimes A) + (A \otimes W_2^\pm)$ are $(H \times H)$-primes of $A \otimes A$, and thus completely prime. Hence, it will suffice to prove that

$$(*) \qquad P_{t,l}^\pm + Q_{t,l'}^\pm = \Delta^{-1}\big((V_1^\pm \otimes A) + (A \otimes W_2^\pm)\big)$$

in all four cases. (There are four cases because $V_1^\pm$ (respectively, $W_2^\pm$) must have the same exponent as $P_{t,l}^\pm$ (respectively, $Q_{t,l'}^\pm$).)

Note that the following collections of $t \times t$ index pairs are hereditary:

$$\{(I,J) \mid |I| = t \text{ and } I \subseteq \{n-l+1, \ldots, n\}\}$$
$$\{(I,J) \mid |I| = t \text{ and } J \subseteq \{1, \ldots, l'\}\}.$$

Hence, Lemma 6.4 shows that

(1) $P_{t,l}^+$ has a basis consisting of those preferred products containing a row $(I, J)$ such that either $|I| > t$, or $|I| = t$ and $I \subseteq \{n-l+1, \ldots, n\}$.
(2) $Q_{t,l'}^-$ has a basis consisting of those preferred products containing a row $(I, J)$ such that either $|I| > t$, or $|I| = t$ and $J \subseteq \{1, \ldots, l'\}$.

Combining (1) and (2) with Lemma 6.5, we obtain

(1') $V_1^+$ has a basis consisting of those preferred products containing a row $(I, J)$ such that either $J \not\subseteq \{1, \ldots, t\}$, or $|I| = t$ and $I \subseteq \{n-l+1, \ldots, n\}$.
(2') $W_2^-$ has a basis consisting of those preferred products containing a row $(I, J)$ such that either $I \not\subseteq \{1, \ldots, t\}$, or $|I| = t$ and $J \subseteq \{1, \ldots, l'\}$.

To get useful bases for the other two cases, we apply the automorphism $\tau$ discussed in (5.4). Observe that $P_{t,l}^- = \tau(Q_{t,l}^-)$ and $Q_{t,l'}^+ = \tau(P_{t,l'}^+)$. Since $\tau$ converts preferred products to antipreferred products, it follows from (1) and (2) that

(3) $P_{t,l}^-$ has a basis consisting of those antipreferred products containing a row $(I, J)$ such that either $|I| > t$, or $|I| = t$ and $I \subseteq \{1, \ldots, l\}$.

(4) $Q_{t,l'}^+$ has a basis consisting of those antipreferred products containing a row $(I, J)$ such that either $|I| > t$, or $|I| = t$ and $J \subseteq \{n - l' + 1, \ldots, n\}$.

Combining (3) and (4) with Lemma 6.5, we obtain

(3') $V_1^-$ has a basis consisting of those antipreferred products containing a row $(I, J)$ such that either $J \nsubseteq \{1, \ldots, t\}$, or $|I| = t$ and $I \subseteq \{1, \ldots, l\}$.

(4') $W_2^+$ has a basis consisting of those antipreferred products containing a row $(I, J)$ such that either $I \nsubseteq \{1, \ldots, t\}$, or $|I| = t$ and $J \subseteq \{n - l' + 1, \ldots, n\}$.

It is clear that $\Delta(P_{t,l}^\pm + Q_{t,l'}^\pm) \subseteq (V_1^\pm \otimes A) + (A \otimes W_2^\pm)$. Hence, if $(*)$ fails, there exists an element

$$x \in \Delta^{-1}\big((V_1^\pm \otimes A) + (A \otimes W_2^\pm)\big) \setminus \big(P_{t,l}^\pm + Q_{t,l'}^\pm\big).$$

Write $x = \sum_{i=1}^m \alpha_i [T_i \mid T_i']$ where the $\alpha_i \in k^\times$ and the $(T_i, T_i')$ are distinct bitableaux of the following types (depending on the four cases):

$(-,-)$: All $(T_i, T_i')$ standard.
$(+,-)$: All $(T_i, T_i')$ preferred.
$(-,+)$: All $(T_i, T_i')$ antipreferred.
$(+,+)$: All $(T_i, T_i')$ antistandard.

We may delete any $[T_i \mid T_i']$ which happen to lie in $P_{t,l}^\pm + Q_{t,l'}^\pm$. Thus, without loss of generality, $[T_i \mid T_i'] \notin P_{t,l}^\pm + Q_{t,l'}^\pm$ for all $i$. In particular, each $T_i$ has at most $t$ columns.

By Proposition 6.3, there must be an index $i$ such that either $[T_i \mid \mu(T_i)] \in V_1^\pm$ or $[\mu(T_i) \mid T_i'] \in W_2^\pm$. Since $T_i$ has at most $t$ columns, all the entries in $\mu(T_i)$ are from $\{1, \ldots, t\}$. If we are in one of the cases $(+, \pm)$, then $T_i$ is decreasing, and so $[T_i \mid \mu(T_i)]$ is a preferred product. Since $V_1^+$ has a basis of preferred products as in (1'), it follows that $[T_i \mid \mu(T_i)] \in V_1^+$ only if $[T_i \mid \mu(T_i)]$ contains a row $(I, J)$ with $|I| = t$ and $I \subseteq \{n - l + 1, \ldots, n\}$, in which case $[T_i \mid T_i'] \in P_{t,l}^+$. On the other hand, in the cases $(-, \pm)$, we have $T_i$ increasing, and so $[T_i \mid \mu(T_i)] = [T_i \mid \mu'(T_i)]$ is an antipreferred product. Then, in view of (3'), we can have $[T_i \mid \mu(T_i)] \in V_1^-$ only if $[T_i \mid T_i'] \in P_{t,l}^-$. To summarize: If $[T_i \mid \mu(T_i)] \in V_1^\pm$, then we would have $[T_i \mid T_i'] \in P_{t,l}^\pm$, which we have ruled out. Similarly, $[\mu(T_i) \mid T_i'] \in W_2^\pm$ would imply $[T_i \mid T_i'] \in Q_{t,l'}^\pm$, another contradiction.

Therefore $(*)$ must hold. $\square$

## 7. $H$-spec $\mathcal{O}_q(M_3(k))$

In this section, we complete the description of the $H$-prime ideals of $\mathcal{O}_q(M_3(k))$, for $q$ not a root of unity. Fix $A = \mathcal{O}_q(M_3(k))$.

**7.1.** In Figure 4.1, we displayed 144 generating sets for ideals of $A$, and in Section 4 we showed that 138 of those ideals are $H$-prime. The cases remaining are the last 6 schematics of Figure 4.2. Four of these schematics generate the ideals $P_{2,2}^{\pm} + Q_{2,2}^{\pm}$ investigated in Theorem 6.6, and these four ideals are $H$-primes.

**7.2.** Let us prove that the ideal

$$P = \langle [I \mid J] \mid I = \{2,3\} \text{ or } J = \{1,2\}\rangle + \langle X_{13}\rangle$$

is an $H$-prime. This will take care of the seventh schematic in Figure 4.2, and the tenth case holds by symmetry (apply the automorphism $\tau$ of (5.4)). Set

$$V_1 = \langle X_{12}, X_{13}, X_{23}, X_{33}, [23 \mid 12]\rangle \quad \text{and} \quad W_2 = \langle X_{13}, X_{31}, X_{32}, X_{33}, [12 \mid 12]\rangle.$$

The factors $A/V_1$ and $A/W_2$ are isomorphic to skew polynomial extensions of the domain $\mathcal{O}_q(M_2(k))/\langle \square \rangle$, so they are domains, and thus $V_1$ and $W_2$ are $H$-primes of $A$. Applying [**6**, (3.1) and Lemma 3.2], we find that $(V_1 \otimes A) + (A \otimes W_2)$ is an $(H \times H)$-prime of $A \otimes A$, and thus completely prime. To show that $P$ is an $H$-prime, it therefore suffices to show that

$$(*) \qquad P = \Delta^{-1}\big((V_1 \otimes A) + (A \otimes W_2)\big)$$

We know from Lemma 6.5 that the ideal $V_0 = \langle X_{13}, X_{23}, X_{33}\rangle$ is spanned by those preferred products $[T \mid T']$ such that 3 occurs in $T'$, that is, $T'$ contains a row $(I, J)$ with $J \not\subseteq \{1, 2\}$. Modulo $V_0$, the elements $[23 \mid 12]$ and $X_{12}$ are normal, and so

$$V_1 = V_0 + [23 \mid 12]A + AX_{12}.$$

Note that $(\{1\}, \{2\})$ (which we shall abbreviate $(1, 2)$) is maximum among index pairs $(I, J)$ with $J \subseteq \{1, 2\}$, and that $(\{2, 3\}, \{1, 2\})$ (which we shall abbreviate $(23, 12)$) is minimum among such index pairs. Hence, if $[T \mid T']$ is any preferred product, either $[T \mid T']X_{12} \in V_0$ or $[T \mid T']X_{12}$ is preferred and similarly, either $[23 \mid 12][T \mid T'] \in V_0$ or $[23 \mid 12][T \mid T']$ is preferred. Therefore $V_1$ is spanned by those preferred products which contain a row $(23, 12)$ or a row $(1, 2)$ or a row $(I, J)$ with $J \not\subseteq \{1, 2\}$.

Similarly, $W_2$ is spanned by those preferred products which contain a row $(12, 12)$ or a row $(1, 3)$ or a row $(I, J)$ with $I \not\subseteq \{1, 2\}$.

It is clear that $\Delta(P) \subseteq (V_1 \otimes A) + (A \otimes W_2)$. Thus if $(*)$ fails, there exists an element $x$ in $\Delta^{-1}\big((V_1 \otimes A) + (A \otimes W_2)\big) \setminus P$. Write $x = \sum_{i=1}^m \alpha_i [T_i \mid T_i']$ where the $\alpha_i \in k^\times$ and the $(T_i, T_i')$ are distinct preferred bitableaux. We may delete any $[T_i \mid T_i']$ which are in $P$. Thus, without loss of generality, $[T_i \mid T_i'] \notin P$ for all $i$.

By case (3) of Proposition 6.3, there is an index $i$ such that

(1) $[T_i \mid \mu(T_i)] \in V_1$ or $[\mu(T_i) \mid T_i'] \in W_2$;
(2) $[T_i \mid \mu_{(12)}(T_i)] \in V_1$ or $[\mu_{(12)}(T_i) \mid T_i'] \in W_2$.

Since $[T_i \mid T_i'] \notin P$, we see that $T_i$ has at most two columns and that $(T_i, T_i')$ contains no row of any of the following forms: $(23, **)$, $(**, 12)$, $(1, 3)$. In particular, all entries of $\mu(T_i)$ and $\mu_{(12)}(T_i)$ are 1's or 2's.

Since $T_i$ does not contain a row $(23)$ and $\mu(T_i)$ does not contain a row $(2)$, we see that $(T_i, \mu(T_i))$ cannot contain either a row $(23, 12)$ or a row $(1, 2)$, and so $[T_i \mid \mu(T_i)] \notin V_1$. Hence, (1) implies that $[\mu(T_i) \mid T_i'] \in W_2$, and so $(\mu(T_i), T_i')$ must contain either a row $(12, 12)$ or a row $(1, 3)$. Since $(T_i, T_i')$ contains no row of the form $(**, 12)$, we see that $(\mu(T_i), T_i')$ cannot contain a row $(12, 12)$. Hence, it must contain a row $(1, 3)$, and thus $(T_i, T_i')$ must have a row of the form $(a, 3)$. Finally, $a \neq 1$ since $(T_i, T_i')$ contains no row $(1, 3)$.

Similarly, $(\mu_{(12)}(T_i), T_i')$ cannot contain either a row $(12, 12)$ or a row $(1, 3)$, whence $[\mu_{(12)}(T_i) \mid T_i'] \notin W_2$, and so $[T_i \mid \mu_{(12)}(T_i)] \in V_1$ by (2). Thus, since $(T_i, \mu_{(12)}(T_i))$ cannot contain a row $(23, 12)$, it must have a row $(1, 2)$, and so $(T_i, T_i')$ contains a row of the form $(1, b)$ where $b \neq 3$.

Since $(T_i, T_i')$ is preferred, either $(a, 3) \leq (1, b)$ or $(1, b) \leq (a, 3)$. However, the first case is impossible since $3 \not\leq_c b$, and the second case is impossible because $1 \not\leq_r a$.

Therefore $(*)$ holds, and the proof is complete.

**7.3.** We have now shown that each schematic in Figure 4.1 generates an $H$-prime of $A$, one which contains the quantum determinant but not all $2 \times 2$ quantum minors. To see that these $H$-primes are distinct, we study their images under suitable retraction maps, as in (3.6).

First, let $\theta_{2,3} : A \to \mathcal{O}_q(M_{2,3}(k))$ be the natural $k$-algebra retraction with kernel $\langle X_{31}, X_{32}, X_{33} \rangle$, where $\mathcal{O}_q(M_{2,3}(k))$ is identified with the subalgebra of $A$ generated the $X_{1j}$ and $X_{2j}$. The images under $\theta_{2,3}$ of the $H$-primes indicated in Figure 4.1 are ideals of $\mathcal{O}_q(M_{2,3}(k))$ which we record in groups in Figure 7.3A. Here the numbers at the left indicate groups of rows in Figure 4.1, while the columns in Figure 7.3A correspond to the columns of Figure 4.1.

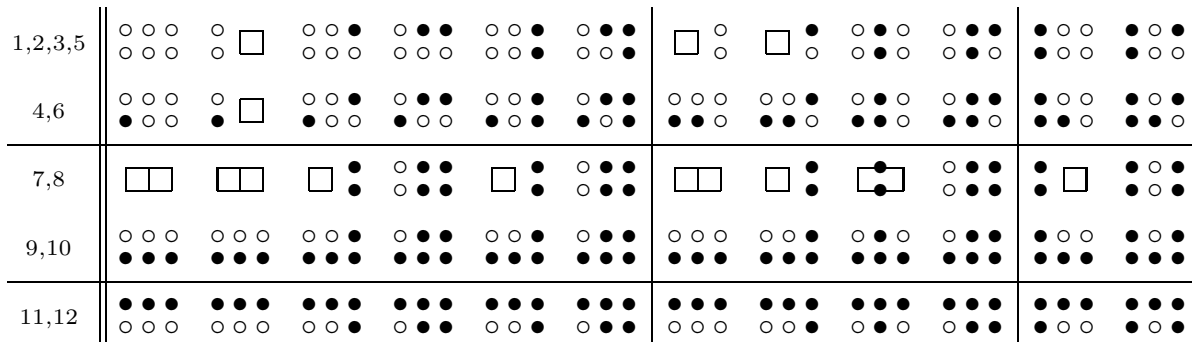

Figure 7.3A

We draw several conclusions from the information in Figure 7.3A:

(a) There are no repetitions from one of the row groups $\{1, 2, 3, 5\}$, $\{4, 6\}$, $\{7, 8\}$, $\{9, 10\}$, $\{11, 12\}$ to another. That is, any $H$-prime arising in one of the rows 1, 2, 3, 5 of Figure 4.1 is distinct from any $H$-prime arising in any of the other rows, and so on.

(b) Within rows 1–6, there are no repetitions from one column to another.

(c) Within columns 7–12, there are no repetitions from one column to another.

By symmetry, we obtain the transposed conditions as well. In particular,

(a') There are no repetitions from one of the column groups $\{1, 2, 3, 5\}$, $\{4, 6\}$, $\{7, 8\}$, $\{9, 10\}$, $\{11, 12\}$ to another.

Thus, if any repetition from one column to another occurs, both columns must lie in the range 1–6, and the repetition must occur within rows 7–12.

Next, let $\theta : A \to k\langle X_{ij} \mid i = 1, 2, 3;\ j = 2, 3\rangle$ be the natural $k$-algebra retraction with kernel $\langle X_{11}, X_{21}, X_{31}\rangle$. We record the images under $\theta$ of the $H$-primes from rows 7–12 and columns 1–6 of Figure 4.1 in Figure 7.3B. From this figure, we see that there cannot be any repetitions from one column to another within the corresponding part of Figure 4.1.

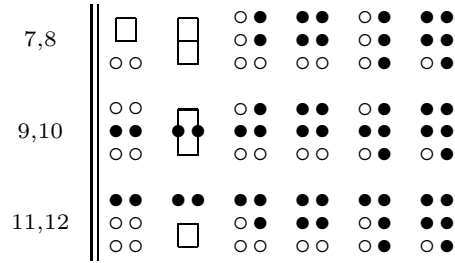

Figure 7.3B

We now conclude that any two $H$-primes of $A$ arising from different columns of Figure 4.1 must be distinct. By symmetry, the same holds for $H$-primes from different rows. Therefore the $H$-primes generated by the 144 schematics in Figure 4.1 are all distinct. Since (by Theorem 1.6) $A$ has exactly 144 $H$-primes which contain the quantum determinant but not all $2 \times 2$ quantum minors, we have found them all.

Our main result is now achieved:

**7.4. Theorem.** *Let $A = \mathcal{O}_q(M_3(k))$ for an arbitrary field $k$ and a non-root of unity $q \in k^\times$, and let $H = (k^\times)^3 \times (k^\times)^3$ act on $A$ by left and right winding automorphisms in the standard manner. Then $A$ has precisely 230 $H$-primes, and they are generated by the schematics displayed in Figures 2.4, 2.5, 3.1, 4.1.* $\square$

We can now conclude that the conjecture discussed in the Introduction holds in the $3 \times 3$ case. Recall that a sequence $u_1, \ldots, u_t$ of elements in a ring is said to be *polynormal* provided $u_i$ is normal modulo the ideal $\langle u_1, \ldots, u_{i-1}\rangle$ for $i = 1, \ldots, t$.

**7.5. Corollary.** *Let $k$ be a field and $q \in k^\times$ a non-root of unity. Then every winding-invariant prime of $\mathcal{O}_q(M_3(k))$ is generated by a set of quantum minors which can be arranged in a polynormal sequence.*

*Proof.* It is clear from Theorem 7.4 that every $H$-prime of $\mathcal{O}_q(M_3(k))$ is generated by a set of quantum minors. In most cases, it is obvious how to arrange the generating set indicated in Figures 2.4, 2.5, 3.1, 4.1 into a polynormal sequence. However, there are a few cases in which some of the 'redundant' quantum minors not drawn in these figures must be included in the generating sets. For example, consider the third schematic in the first row of Figure 2.5, which indicates 6 quantum minors (three $1 \times 1$'s and three $2 \times 2$'s). None of the $X_{i2}$ is normal modulo the ideal generated by the three indicated $2 \times 2$ quantum minors, but the $X_{i2}$ are normal modulo the ideal generated by all the $2 \times 2$ quantum minors. Thus one obtains a polynormal sequence of 12 quantum minors generating this $H$-prime, beginning with all nine $2 \times 2$ quantum minors in a suitable order. We leave the remaining details of the proof to the reader. □

Department of Mathematics, University of California, Santa Barbara, CA 93106, USA
*E-mail address*: goodearl@math.ucsb.edu

Department of Mathematics, J.C.M.B., Kings Buildings, Mayfield Road, Edinburgh EH9 3JZ, Scotland
*E-mail address*: tom@maths.ed.ac.uk